\newif\ifjournal
\journaltrue
\newif\ifarxiv
\arxivtrue
\newif\ifsinglcol
\singlcoltrue

\ifsinglcol
	\ifarxiv
		\documentclass[10pt,onecolumn,twoside,draftclsnofoot]{IEEEtran}
	\else
		\documentclass[10pt,draftcls,twoside,onecolumn]{IEEEtran}
	\fi
\else
\documentclass[journal]{IEEEtran}
\fi

\makeatletter

\usepackage{amsthm}
\usepackage{enumitem}

\usepackage{graphicx}      
\graphicspath{./Matlab/}	
\usepackage[usenames,dvipsnames]{color}

\usepackage{mathtools}
\usepackage{cite}
\usepackage{amssymb}
\usepackage[mathscr]{eucal}
\usepackage{units}
\usepackage{yhmath}
\usepackage{units}
\usepackage{pgfplots}
\usepackage[utf8]{inputenc}
\usepackage{balance}
\usepackage{dsfont}
\usepackage{caption}
\usepackage{subcaption}
\usepackage{dblfloatfix}
\usepackage{color} 

\newcommand{\added}[1]{#1}

\usepackage{tikz}
\usetikzlibrary{decorations.pathmorphing} 
\usetikzlibrary{matrix} 
\usetikzlibrary{arrows} 
\usetikzlibrary{calc} 

\tikzstyle{block} = [draw,rectangle,thick,minimum height=2em,minimum width=2em]
\tikzstyle{sum} = [draw,circle,inner sep=0mm,minimum size=2mm]
\tikzstyle{connector} = [->,thick]
\tikzstyle{line} = [thick]
\tikzstyle{branch} = [circle,inner sep=0pt,minimum size=1mm,fill=black,draw=black]
\tikzstyle{guide} = []
\tikzstyle{snakeline} = [connector, decorate, decoration={pre length=0.2cm,
                         post length=0.2cm, snake, amplitude=.4mm,
                         segment length=2mm},thick, magenta, ->]

\def\myarm{1cm}
\def\myangle{0}
\tikzset{
  arm/.default=1cm,
  arm/.code={\def\myarm{#1}}, 
  angle/.default=0,
  angle/.code={\def\myangle{#1}} 
}

\usepackage{blkarray}
\newcommand{\matindex}[1]{\mbox{\scriptsize $#1$}}

\newtheorem{theorem}{Theorem}
\newtheorem{assumption}{Assumption}

\newtheorem{lemma}[theorem]{Lemma}
\newtheorem{remark}[theorem]{Remark}
\newtheorem{problem}{Problem}

\DeclareMathOperator{\diag}{diag}
\DeclareMathOperator{\blockdiag}{blockdiag}

\DeclareMathOperator*{\argmax}{arg\,max}
\newcommand{\1}{\ensuremath{\mathbf{1}}}
\newcommand{\0}{\ensuremath{\mathbf{0}}}

\newcommand{\G}{\ensuremath{\mathcal{G}}}
\renewcommand{\L}{\ensuremath{\mathcal{L}}}
\newcommand{\Z}{\ensuremath{\mathcal{Z}}}
\newcommand{\B}{\ensuremath{\mathcal{B}}}
\newcommand{\N}{\ensuremath{\mathcal{N}}}
\newcommand{\E}{\ensuremath{\mathcal{E}}}
\newcommand{\K}{\ensuremath{\mathcal{K}}}

\newcommand{\R}{\ensuremath{\mathds{R}}}
\newcommand{\D}{\ensuremath{\mathcal{D}}}

\newcommand{\ddx}[1]{\frac{\partial}{\partial #1}}
\newcommand{\dddx}[1]{\frac{\partial^2}{\partial #1^2}}
\newcommand{\dd}[2]{\frac{\partial #1}{\partial #2}}
\allowdisplaybreaks

\newboolean{showcomments}
\setboolean{showcomments}{true}
\newcommand{\slow}[1]{\ifthenelse{\boolean{showcomments}}
{ \textcolor{red}{(Steven says:  #1)}}{}}
\newcommand{\emallada}[1]{\ifthenelse{\boolean{showcomments}}
{ \textcolor{red}{(Enrique says:  #1)}}{}}

\newcommand{\hide}[1]{}

\title{Optimal load-side control for frequency regulation in smart grids}

\ifarxiv
\author{Enrique Mallada, Changhong Zhao, and Steven Low
\thanks{
This work was supported by NSF NetSE grant CNS 0911041, ARPA-E grant DE-AR0000226, Southern California Edison, National Science Council of Taiwan R.O.C. grant NSC 103-3113-P-008-001, and Caltech Resnick Institute.}
\thanks{The authors are with the Department of Computing + Mathematical Sciences, California Institute of Technology, Pasadena, CA 91125 USA (e-mail: \{mallada,czhao,slow\}@caltech.edu).}
\thanks{}
}
\else

\author{Enrique Mallada,~\IEEEmembership{Member,~IEEE,} Changhong Zhao,~\IEEEmembership{Student Member,~IEEE,} and Steven Low,~\IEEEmembership{Fellow,~IEEE}
\ifjournal
\thanks{{\bf Reserved space.-}}
\thanks{}\thanks{}\thanks{}
\fi
\thanks{
This work was supported by NSF NetSE grant CNS 0911041, ARPA-E grant DE-AR0000226, Southern California Edison, National Science Council of Taiwan R.O.C. grant NSC 103-3113-P-008-001, and Caltech Resnick Institute.}
\thanks{The authors are with the Department of Computing + Mathematical Sciences, California Institute of Technology, Pasadena, CA 91125 USA (e-mail: \{mallada,czhao,slow\}@caltech.edu).}
\thanks{}
}
\fi

\ifjournal
\ifsinglcol\else
\renewcommand{\baselinestretch}{1}
\fi
\else
\renewcommand{\baselinestretch}{.955}
\fi

\begin{document}

\maketitle
\begin{abstract}                
Frequency control rebalances supply and demand while maintaining the network state within operational margins. It is implemented using fast ramping reserves that are expensive and wasteful, and which are expected to grow with the increasing penetration of renewables. The most promising solution to this problem is the use of demand  response, i.e. load  participation in frequency control. Yet it is still unclear how to efficiently integrate load participation without introducing instabilities and violating operational constraints.

In this paper we present a comprehensive load-side frequency control mechanism that can maintain the grid within operational constraints. In particular, our controllers can rebalance supply and demand after disturbances, restore the frequency to its nominal value and preserve inter-area power flows. Furthermore, our controllers are distributed (unlike the currently implemented frequency control), can allocate load updates optimally, and can maintain line flows within thermal limits.
We prove that such a distributed load-side control is globally asymptotically stable and robust to unknown load parameters. We illustrate its effectiveness through simulations.
\end{abstract}



\section{Introduction}
\label{sec:introduction}

Frequency control maintains the frequency of a power network at its nominal
value when demand or supply fluctuates.   It is traditionally implemented on the generation side
and consists of three mechanisms that work in concert \cite{WoodWollenberg1996, Bergen2000, MachowskiBialek2008}.
The primary frequency control, called the droop control and completely decentralized,
operates on a timescale up to low tens of 
seconds and uses a governor to adjust, around a setpoint, the mechanical power input to a generator 
based on the local frequency deviation.   
The primary control can rebalance power and stabilize the frequency    
but does not restore the nominal frequency.   The secondary frequency control
(called automatic generation control or AGC) operates on a timescale up to a minute or so and adjusts the
setpoints of governors  in a control area  in a centralized fashion to drive 
the frequency back to its nominal value and the inter-area power flows to their scheduled values.  
Finally, economic dispatch operates on a timescale of several minutes or up and schedules the output levels 
of generators that are online and the inter-area power flows.
See \cite{ilic2007hierarchical, kiani2012hierarchical} for a recent hierarchical model of power systems
and their stability analysis.

Load-side participation in frequency control offers many advantages, including faster response, lower
fuel consumption and emission, and better localization of disturbances.  The idea of using
frequency adaptive loads dates back to \cite{Schweppe1980} that 
advocates their large scale deployment to 
``assist or even replace turbine-governed systems and spinning reserve.''   
They also proposed to use spot prices to incentivize the users to adapt their consumption 
to the true cost of generation at the time of consumption.   
Remarkably it was emphasized back then that such frequency adaptive loads will
``allow the system to accept more readily a stochastically fluctuating energy source, such as wind or solar 
generation'' \cite{Schweppe1980}.

This is echoed recently in, e.g., 
\cite{Trudnowski2006, LuHammerstrom2006,Short2007, donnelly2010frequency, brooks2010demand, Callaway:2011gg,molina2011decentralized} 
that argue for ``grid-friendly'' appliances, such as refrigerators, water or space heaters, ventilation systems, and air 
conditioners, as well as plug-in electric vehicles to help manage energy imbalance.
Simulations in all these studies have consistently shown significant improvement in performance and
reduction in the need for spinning reserves.
A small scale project by the Pacific Northwest National Lab in 2006--2007 demonstrated the use of
200 residential appliances  in primary frequency control that
automatically reduce their consumption when the frequency of the household dropped below a threshold (59.95Hz) 
\cite{Hammerstrom2007}. \added{Although these simulation studies and field trials are insightful, they fall short in predicting the (potential) behavior of a large-scale deployment of frequency control.}

\added{This has motivated the recent development of new analytic studies assessing the effect of distributed frequency control in power systems~\cite{Andreasson:2014jv,Andreasson:2013ve,Zhang:2013wy,Zhang:2015kq,Zhao:2014bp,Li:2014we,mallada2014distributed,you2014reverse}, and microgrids~\cite{Liu:2014gg,Shafiee:2014eb,Dorfler:2014uu,SimpsonPorco:2013dn,Dorfler:uo}, which can be grouped into three main approaches.
The first approach builds on \emph{consensus algorithms} to provide efficiency guarantees to classical PI controllers~\cite{Andreasson:2014jv,Andreasson:2013ve,Shafiee:2014eb,Dorfler:2014uu,SimpsonPorco:2013dn,Dorfler:uo}. It achieves efficiency but does not manage congestion, i.e., it does 
not enforce constraints such as thermal limits.
The second approach \emph{reverse engineers} the network dynamics as a primal-dual algorithm of an underlying optimization problem, and then add constraints and modifies the objective function while preserving the primal-dual interpretation of the network dynamics~\cite{Zhao:2014bp,Zhao:2013ja,Li:2014we,mallada2014distributed,you2014reverse}. 
It successfully achieves efficiency but it limits the type of operational constraints that can be satisfied.
The third approach directly formulates an optimization problem that encodes \emph{all operational constraints} and then designs a \emph{hybrid algorithm} that combines the network dynamics with a subset of the primal-dual algorithm~\cite{Zhang:2013wy,Zhang:2015kq}. It is able to satisfy operational constraints, but the stability depends on the network parameters.
}
\color{black}
{\it Contributions of this work:}
In this paper we develop a method to achieve secondary frequency regulation
and congestion management (bringing line flows to within their limits), \emph{in a
distributed manner using controllable loads}. To contrast with the
generation-side AGC frequency control, we refer to our solution as \emph{load-side control}.    
To our knowledge, this method produces the \emph{first} distributed controllers for 
demand response that are scalable and enforce required operational 
constraints for frequency regulation, such as restoring nominal frequency and 
preserving inter-area flows, while respecting line limits. 
\color{black}

\added{
Our work builds on previous optimization-based approaches~\cite{Zhao:2014bp,Zhao:2013ja,Li:2014we, mallada2014distributed, you2014reverse, Zhang:2013wy, Zhang:2015kq}.
The crux of our solution 
is the introduction of \emph{virtual (line) flows} that can be used to implicitly 
constrain real flows without altering the primal-dual interpretation of the 
network dynamics.
A virtual flow is a cyber quantity on each line that a controller computes 
	based on information from its neighbors.  In steady state, its value equals the 
	actual line flows incident on that controller.  This device allows us to impose 
	arbitrary constraints on the (actual) line flows for congestion management
	and restoring inter-area flows, while retaining the ability to
	exploit network dynamics to help carry out the primal-dual algorithm. 
}


Our contribution with respect to the existing literature is manifold.
Unlike \cite{Andreasson:2013ve,Zhang:2013wy,Zhang:2015kq}, our global asymptotic stability result (Theorem \ref{thm:global-convergence} in Section \ref{sec:analysis}) is independent of the controller gains, which is highly desirable for fully distributed deployments. Our results hold for arbitrary topologies, in contrast with~\cite{Zhang:2013wy,Dorfler:2014uu}, and can impose inter-area constraints, thermal limits~\added{or \emph{any linear equality or inequality constraint} in the line flows}. Moreover, we provide a convergence analysis in the presence of unknown parameters (Section \ref{sec:implementation-details}) that is novel in the primal-dual literature and provides the necessary robustness for large scale distributed deployments.
\added{Finally, our framework can further extend to include intermediate buses without generators or loads -- quite common in transmission networks -- which are not considered by the existing literature, and to fully distribute non-local constraints like those imposed on inter-area flows (Section \ref{sec:framework-extensions}).}

A preliminary version of this work has been presented in~\cite{Mallada:2014ej}. This paper extends~\cite{Mallada:2014ej} in several ways. 
First, the robustness study of our controllers with respect to uncertain load parameter (Section \ref{sec:implementation-details}) as well as the framework extensions (Section \ref{sec:framework-extensions}) are new.
Second, we include detailed proofs that were omitted in~\cite{Mallada:2014ej} due to space constraints. 
Finally, we extend our simulations in Section \ref{sec:numerical-illustrations} to further
 illustrate the conservativeness of the uncertainty bounds of Theorem \ref{thm:global-convergence-2}.






%
%


\section{Preliminaries}
\label{sec:preliminaries}

Let $\R$ be the set of real numbers and $\mathds N$ the set of natural numbers. Given a finite set $S\subset \mathds N$ we use $|S|$ to denote its cardinality. For a set of
scalar numbers $a_i\in \R$, $i\in S$, we denote $a_S$ to its column vector, i.e. $a_S:=(a_i)_{i\in S}\in\R^{|S|}$; we usually drop the subscript $S$ when the set is known from the context. For two vectors $a\in\R^{|S|}$ and $b\in\R^{|S'|}$ we define the column vector $x=(a,b)\in\R^{|S|+|S'|}$.
Given a matrix $A$, we denote its transpose as $A^T$ and use $A_i$ to denote the $i$th row of $A$. We will also use $A_S$ to denote the sub matrix of $A$ composed only of the rows $A_i$ with $i\in S$. 
 The diagonal matrix of a sequence $\{a_i,\;i\in S\}$ is represented by $\diag(a_i)_{i\in S}$. Similarly, for a sequence of matrices $\{A_h,\;h\in S\}$ we let $\blockdiag(A_h)_{h\in S}$ denote the block diagonal matrix. 
Finally, we use either $\1_n$ and $\1_{n\times m}$ ($\0_n$ and $\0_{n\times m}$) to denote the vector and matrix of all ones (zeros), or $1$ ($0$) when its dimension is understood from the context.

\added{
For a function $f:\R^n\rightarrow \R^n$ we use $f'(x):=\ddx{x}f(x)$ to denote the Jacobian and $f^{-1}$ to denote its inverse. When $n=1$, $f''(x)$ denotes the second derivative $\frac{\partial^2}{\partial x^2}f(x)$.}

For given vectors $u\in\mathds{R}^n$ and $a\in\mathds{R}^n$,  and set $S\subseteq\{1,\dots,n\}$, the operator $[a]^+_{u_S}$ is defined element-wise by  
\begin{equation}\label{eq:projection}
([a]_{u_S}^+)_i=\begin{cases}
[a_i]_{u_i}^+, & \mbox{ if } i\in S,\\
a_i, & \mbox{ if } i\not\in S,
\end{cases}
\end{equation}
where  $[a_i]_{u_i}^+$ is equal to $a_i$ if either $a_i>0$ or $u_i>0$, and $0$ otherwise.
Whenever  $u_S^*\geq 0$, the following relation holds:
\begin{equation}\label{eq:projection-property}
(u_S-u_S^*)^T[a_S]_{u_S}^+\leq (u_S-u^*_S)^Ta_S
\end{equation}
since for any pair $(u_i,a_i)$, with $i\in S$, that makes the projection active ($[a_i]_{u_i}^+=0$) we must have $u_i\leq0$ and $a_i\leq0$, and therefore 
$
(u_i-u_i^*)a_i\geq 0 = (u_i-u_i^*)^T[a_i]_{u_S}^+.
$

\subsection{Network Model}\label{ssec:network-model}

We consider a power network described by a directed graph $G(\N,\E)$ where $\N=\{1,...,|\N|\}$ is the set of buses, denoted by either $i$ or $j$, and $\E\subset \N\times\N$ is the set of transmission lines denoted by either $e$ or $ij$ such that if $ij\in \E$, then $ji\not\in \E$. 

We partition the buses $\N=\G \cup \L$ and use $\G$ and $\L$  to denote the set of generator and load buses respectively. We assume that a generator bus may also have loads attached to it, but not otherwise.
We assume the graph $G(\N,\mathcal{E})$ is connected, and make the following assumptions which are standard and well-justified for transmission networks \cite{Kundur:1994tx}:
(i) Bus voltage magnitudes $|V_i|=|V_j|=1\text{pu}$ for $i,j \in \N$. 
(ii) Lines $ij \in \mathcal E$ are lossless and characterized by their susceptances $B_{ij}>0$. 
(iii) Reactive power flows do not affect bus voltage phase angles and frequencies.

The evolution of the transmission network is described by 
\begin{subequations}\label{eq:swing-theta}
\begin{align}
\dot \theta_i &\!=\! \omega_i  & i\in\N\label{eq:sw-theta-c} \\
M_i \dot\omega_i & \!=\! P^{in}_i \!\!-\!(d_i \!+\! D_i\omega_i) \!-\!\textstyle{\sum_{j\in \N_i}\! B_{ij}(\theta_i\!-\!\theta_j)}& i \in \G \label{eq:sw-theta-a}\\
0 & \!=\! P^{in}_i\!\!-\!(d_i \!+\! D_i\omega_i) \!-\!\!\textstyle\sum_{j\in \N_i}\! B_{ij}(\theta_i\!-\!\theta_j)& i \in \L \label{eq:sw-theta-b}
\end{align}
\end{subequations}
where  $d_i$  denotes an aggregate controllable load, $D_i\omega_i$ denotes the aggregate frequency sensitive consumption
due to uncontrollable but frequency-sensitive loads and/or generators' damping loss. 
$M_i$ denotes the generator's inertia, \added{and $P^{in}_i$ is the difference between mechanical power injected by a generator and the constant aggregate power consumed by loads. Notice that for load buses ($i\in \L$) $P^{in}_i<0$.}

\added{
Equation \eqref{eq:swing-theta}  represents the standard swing equations when the bus frequency $\omega_i$ and phase $\theta_i$ are close to nominal values $\omega^0$ and $\theta_{i}^0$.}\footnote{Without loss of generality, we take here $\omega_0=0$.} 
\added{
Given \eqref{eq:swing-theta}, however, our results hold for arbitrary $\omega_i$ and $\theta_i$. 
}
By letting $P_e=P_{ij}:=B_{ij}(\theta_i-\theta_j)$ for $e=ij\in\E$, we can equivalently rewrite \eqref{eq:swing-theta} as
\begin{subequations}\label{eq:1}
\begin{align}
M_i \dot\omega_i &  \!=\! P^{in}_i -(d_i + D_i\omega_i) -\textstyle\sum_{e\in \E} C_{i,e}P_{e} \;& i \in \G \label{eq:1a}\\
0 & \!=\!P^{in}_i -(d_i + D_i\omega_i) -\textstyle\sum_{e\in \E} C_{i,e}P_{e} & i \in \L \label{eq:1b}\\
\dot P_{ij} &\!=\!  B_{ij}(\omega_i - \omega_j) & ij\in\E\label{eq:1c}
\end{align}
\end{subequations}
where $C_{i,e}$ are the elements of the incidence matrix $C\in \R^{|\N|\times|\E|}$ of the graph $G(\N,\E)$ such that $C_{i,e}=-1$ if $e=ji\in\E$, $C_{i,e}=1$ if $e=ij\in \E$ and $C_{i,e}=0$ otherwise, and the line flows initial condition must satisfy $P_{ij}(0)=B_{ij}(\theta_i^0-\theta_j^0)$.

\begin{remark}\label{rem:zero-injection}
Our model assumes that every bus has both $D_i\!>\!0$ and controllable load $d_i$. While this is reasonable for load and generator buses (for $i\in\G$ the generator can implement $-d_i$), it is unreasonable for intermediate buses that have neither generators nor loads. This is addressed in Section \ref{sec:framework-extensions}. 
Our framework can also handle the case where $d_i$ is present only for $i\in S\subset\N$. This case is omitted due to space constraints.
\end{remark}

\hide{
For notational convenience we will use whenever needed the vector form of \eqref{eq:1}, i.e.
\begin{align*}
M_\G \dot\omega_\G & = P^{in}_\G -(d_\G + D_\G\omega_\G) -C_\G P\\
0 & = P^{in}_\L-(d_\L + D_\L\omega_\L) -C_\L P\\
\dot P &=  D_BC^T\omega
\end{align*}
where the matrices $C_\L$ and $C_\G$ are defined by splitting the the rows of $C$ between generator and load buses, i.e. $C=[C_\G^T\;C_\L^T]^T$, $D_B=\diag(B_{ij})_{ij\in\E}$, $M_\G=\diag(M_i)_{i\in\G}$ and for a given set $S\subset \N$ $D_S=\diag(D_i)_{i\in S}$.
}

\hide{
Equation \eqref{eq:1} is a non-standard formulation of the swing equations~\addcites~ where in stead of phase angles $\theta_i$ we use power flows $P_{ij}=B_{ij}(\theta_i-\theta_j)$ as state variables. This change of variable allows us to interpret the swing equations as an primal-dual algorithm solving an underlying optimization problem~\cite{Zhao:2013ja}. We refer de reader to Appendix \ref{app:swing-equations} for a detailed explanation of the connection between \eqref{eq:1} and the phase model and its extension for networks with constant resistance-reactance ratios.
}

\subsection{Operational Constraints}\label{ssec:operational-constraints}

We denote each control area using $k$ and let $\K:=\{1,\dots,|\K|\}$ be the set of all areas in the network.
Let $\N_k\subseteq\N$ be set of buses in area $k$ and $\B_k\subseteq \E$ the set of boundary edges of area $k$, i.e. $ij\in\B_k$ iff either $i\in\N_k$ or $j\in\N_k$, but not both.

Within each area, the AGC scheme seeks to restore the frequency to its nominal value and preserve a constant power transfer outside the area, i.e.
\begin{equation}\label{eq:inter-area-constraint-2}
\textstyle\sum_{ij\in\B_k} \hat C_{k,ij} P_{ij} =\hat P_k
\end{equation}
where $\hat C_{k,ij}$ is equal to $1$, if $i\in \N_k$, $-1$, if $j\in \N_k$, and $\hat P_k$ is the net scheduled power injection of area $k$.

By defining $\hat C_{k,ij}$ to be $0$ whenever $ij\not\in\B_k$ we can also relate $\hat C\in\R^{|\K|\times|\E|}$ with $C\in \R^{|\N|\times |\E|}$ using 
\begin{equation}\label{eq:Cflat}
\hat C=E_\K C
\end{equation}
where $E_\K :=[e_1\;\dots\; e_{|\K|}]^T$ and $e_k\in\R^{|\N|}$, $k\in\K$, is a vector with elements $(e_k)_i=1$ if $i\in\N_k$ and $(e_k)_i=0$ otherwise.


Finally, the thermal limit constraints are usually given by
\begin{equation}\label{eq:thermal-limit}
\underline P \leq P\leq \overline P
\end{equation}
where $P:=(P_e)_{e\in\E}$, $\overline P:=(\overline P_e)_{e\in\E}$ and $\underline P:=(\underline P_e)_{e\in\E}$ represent the line flow limits; i.e. $\underline P=-\overline P$ so that $|P|\leq \overline P$.

\section{Optimal load-side control}
\label{sec:frequency-preserving-olc}

Suppose the system \eqref{eq:1} is in equilibrium, i.e. $\dot{\omega}_i=0$ for all $i$ and $\dot{P}_{ij}=0$ for all $ij$,  and at time 0, there is a disturbance, represented by a step change in the vector $P^{in} := (P_i^{in})_{i\in \N}$, that produces a power imbalance.
Then, we are interested in designing a distributed control mechanism that rebalances the system, restores the frequency to its nominal value while maintaining the operational constraints of Section \ref{ssec:operational-constraints}.
Furthermore, we would like this mechanism to produce an efficient allocation among all the users (or loads) that are willing to adapt.

\added{
We use $c_i(d_i)$ to denote the cost or disutility of changing the load consumption by $d_i$. This allows us to formally describe our notion of efficiency in terms of the loads' power share. More precisely, we shall say that a load control is efficient if in equilibrium solves the Optimal Load Control (OLC) problem:}
\textcolor{black}{
\begin{problem}[OLC]\label{problem:olc}
\begin{subequations}\label{eq:welfare}
\begin{flalign}
\underset{d,\omega,\theta}{\text{\rm minimize}} \qquad\qquad\qquad \sum_{i\in\N} &c_i(d_i)   \label{eq:welfare-a}\\
\intertext{\rm subject to}
P^{in} - (d + D\omega) &= L_B\theta&\label{eq:welfare-b}\\
\omega &= 0&\label{eq:welfare-c}\\
\hat C BC^T\theta &= \hat P&\label{eq:welfare-d}\\
\underline P\leq BC^T\theta &\leq \overline P&\label{eq:welfare-e}
\end{flalign}
\end{subequations}
\end{problem}\noindent
where $d=(d_i)_{i\in\N}$, $\omega=(\omega_i)_{i\in\N}$, $\theta=(\theta_i)_{i\in\N}$,  $D=\diag(D_i)_{i\in\N}$, $B=\diag(B_{ij})_{ij\in\E}$, $(BC^T\theta)_{ij}=B_{ij}(\theta_i-\theta_j)$ and $L_B:=CBC^T$ is the $B_{ij}$-weighted Laplacian matrix.
}

Throughout this paper we make the following assumptions:
\begin{assumption}[Cost function]\label{as:1}
The cost function $c_i(d_i)$ is $\alpha$-strongly convex and second order continuously differentiable ($c_i\in C^2$ with $c_i''(d_i)\geq \alpha>0$) in the interior of its domain $\D_i:=[\underline{d}_i,\overline{d}_i]\subseteq \R$, such that $c_i(d_i)\rightarrow +\infty$ whenever $d_i\rightarrow \partial \D_i$.
\end{assumption}

\begin{assumption}[Strict feasibility in $\D$]\label{as:2}
The OLC problem \eqref{eq:olfc} is feasible and there is at least one feasible $(d,\omega,\theta)$ such that $d\in \mathrm{Int} \mathcal\, \D:=\Pi_{i=1}^{|\N|}\mathrm{Int}\, \D_i$.
\end{assumption}
\added{
Assumption \ref{as:2} guarantees that, even in the presence of Assumption \ref{as:1}, the optimal solution of OLC is finite, and therefore allows us to use the KKT conditions~\cite[Ch.\ 5.2.3]{Boyd:2004cv} to characterize it.
}

\added{
Problem \ref{problem:olc} is convex and therefore can be efficiently solved using several optimization algorithms. 
However, unlike standard optimization problems, we can only modify the loads $d_i$ while $\theta_{i}$ and $\omega_i$ react to these changes according to \eqref{eq:swing-theta}.
We overcome these restrictions by formulating an equivalent optimization problem whose primal-dual optimization algorithm embeds the line flow version of the swing equations~\eqref{eq:1}.}

%

\subsection{Virtual Flows Reformulation}
\added{ We now proceed to describe the optimization problem that will be used to derive the distributed controllers that achieve our goals. The crux of our solution comes from implicitly imposing the constraints \eqref{eq:welfare-c}-\eqref{eq:welfare-e} by using \emph{virtual flows} instead of explicitly using $\omega_i$ and $\theta_i$. This, together with an additional quadratic objective on $\omega_i$ and substituting $B_{ij}(\theta_i-\theta_j)$ with $P_e$ in \eqref{eq:welfare-b}, allows us to embed the network dynamics as part of the primal-dual algorithm while preserving all the desired constraints.}
\noindent
\begin{problem}[VF-OLC]\label{problem:olc-aux}
\begin{subequations}\label{eq:olfc}
\begin{flalign}
\underset{ d,\omega,\phi,P}{\text{\rm minimize}} \qquad\quad\;\;\; \sum_{i\in\N} c_i(d_i) +& \frac{D_i{\omega_i}^2}{2} \label{eq:olfc-a}\\
\intertext{\rm subject to}
P^{in} - (d + D\omega) &= CP&\label{eq:olfc-b}\\
P^{in} -d & =L_B\phi&\label{eq:olfc-c}\\
\hat C BC^T\phi &= \hat P&\label{eq:olfc-d}\\
\underline P\leq BC^T\phi&\leq \overline P&\label{eq:olfc-e}
\end{flalign}
\end{subequations}
\end{problem}\noindent
where $\phi=(\phi_i)_{i\in\N}$ represents the virtual phases and $(BC^T\phi)_{ij}=B_{ij}(\phi_i-\phi_j)$ is the corresponding \emph{virtual flow} through line $ij\in\E$.

Although not clear at first sight, the constraint  \eqref{eq:olfc-c} implicitly enforces that any optimal solution of VF-OLC $(d^*,\omega^*,\phi^*,P^*)$ will restore the frequency to its nominal value, i.e. $\omega_i^* = 0$. Thus the additional term on the objective function does not influence the optimal solution. Similarly, we will use constraint \eqref{eq:olfc-d} to impose \eqref{eq:welfare-d}  and \eqref{eq:olfc-e} to impose~\eqref{eq:welfare-e}. 

We use $\nu_i$, $\lambda_i$ and $\pi_{k}$ to denote the Lagrange multipliers of constraints \eqref{eq:olfc-b}, \eqref{eq:olfc-c} and  \eqref{eq:olfc-d},  and $\rho_{ij}^+$ and $\rho_{ij}^-$ as multipliers of the right and left constraints of  \eqref{eq:olfc-e}, respectively.
In order to make the presentation more  compact sometimes we will use $x=(\phi,P)\in\R^{|\N|+|\E|}$ and $\sigma=(\lambda,\nu,\pi,\rho^+,\rho^-)\in\R^{2|\N|+|\K|+2|\E|}$, as well as $\rho:=(\rho^+,\rho^-)$.

Using this notation we can write the Lagrangian of VF-OLC as
\begin{align}
&L(d,\omega,x,\sigma)= \sum_{i\in\N} \left(c_i(d_i) + \frac{D_i{\omega_i}^2}{2}\right)+ \nu^T(P^{in} - (d + D\omega)  \nonumber\\
& \quad-CP)+  \lambda^T(P^{in}-d -L_B\phi)+ \pi^T(\hat C BC^T\phi - \hat P )\nonumber\\
& \quad  + \rho^{+T}(BC^T\phi-\overline P)+ \rho^{-T}(\underline P - BC^T\phi)\allowbreak\nonumber\\
&=\sum_{i\in\N}  c_i(d_i) \!-\! (\lambda_i+\nu_i)d_i \!+\! D_i\omega_i\left({{\omega}_i}/{2}-\nu_i\right)\!+\!(\nu_i+\lambda_i)P_i^{in} \nonumber\\
&\quad-P^TC^T\nu -\phi^T(L_B\lambda -CB\hat C^{T}\pi- CB(\rho^+-\rho^-)) \nonumber\\
&\quad-\pi^T\hat P-\rho^{+T}\overline P + \rho^{-T}\underline P \label{eq:L}
\end{align}

The next lemmas characterize the optimality conditions of VF-OLC and its equivalence with OLC. Their proofs can be found in the Appendix.

\begin{lemma}[Optimality]\label{lem:OLC-optimality}
Let $G(\N,\E)$ be a connected graph. Then $(d^*,\omega^*,\phi^*,P^*,\sigma^*)$ is a primal-dual optimal solution to VF-OLC
if and only if $(d^*,\omega^*,\phi^*,P^*)$ is primal feasible, $\rho^{+*},\rho^{-*}\geq 0$,
\begin{gather}
 \omega_i^* = \nu_i^*,\;\;d_i^*=c'^{-1}_i(\nu_i^*+\lambda_i^*),\;\;\nu_i^*=\hat\nu,\;\; i\in\N, \label{eq:16}
\end{gather}
where $c'^{-1}_i$ is the inverse of the derivative of $c_i$, $\hat \nu$ is some scalar, 
\begin{align}
\textstyle\sum_{j\in\N_i} B_{ij}(\lambda_j^*-\lambda_i^*)+ C_{i}B(\hat C^T\pi^* + \rho^{+*}-\rho^{-*})=0\label{eq:16b}
\end{align}
with $C_i$ being the $i$th row of $C$, and
\begin{subequations}\label{eq:complementary-slackness}
\begin{align}
\rho_{ij}^{+*}(B_{ij}(\phi_i^*-\phi_j^*) - \overline P_{ij}) = 0, \quad ij\in\E,\\ 
\rho_{ij}^{-*}(\underline P_{ij}  - B_{ij}(\phi_i^*-\phi_j^*) ) = 0, \quad ij\in\E.
\end{align}
\end{subequations}
Moreover, $d^*$, $\omega^*$, $\nu^*$ and $\lambda^*$ are unique with $\omega_i^*=\nu_i^*=\hat\nu=0$.
\end{lemma}

\begin{lemma}[OLC and VF-OLC Equivalence]\label{lem:equivalence}
Given any set of vectors $(d^*,\omega^*,\theta^*,\phi^*,P^*)$ satisfying $C^T\theta^*=C^T\phi^*$ and $L_B\theta^* = CP^*$. Then $(d^*,\omega^*,\theta^*)$ is an optimal solution of OLC if and only $(d^*,\omega^*,\phi^*,P^*)$ is an optimal solution to VF-OLC. 
\end{lemma}

\subsection{Distributed Optimal Load-side Control}
\label{ssec:distributed-optimal-load-control}

\color{black}
We now show how to leverage the power network dynamics to solve the OLC problem in a distributed way. 
Our solution is based on the classical primal dual optimization algorithm that has been of great use to design congestion control mechanisms in communication networks~\cite{kelly1998rate,low1999optimization,srikant2004mathematics,palomar2006tutorial}.

Since by Lemma \ref{lem:equivalence}, VF-OLC provides the same optimal load schedule as OLC, we can solve  VF-OLC instead. This will allow us to incorporate the network dynamics as part of an optimization algorithm that indirectly solves OLC.

To achieve this goal we first minimize \eqref{eq:L} over $d$ and $\omega$ which is achieved by setting $\omega_i=\nu_i$ and $d_i=c_i^{'-1}(\nu_i+\lambda_i)$ in \eqref{eq:L}. Thus we get  
\begin{align}
&L(x,\sigma) =\underset{d, \omega}{\text{minimize}} \quad L(d,\omega, x, \sigma) \nonumber\\
&=\Phi_{\N}(\lambda,\nu) - P^TC^T\nu   - \pi^T\hat P - \rho^{+T}\overline{P} + \rho^{-T}\underline{P}\nonumber\\
& \quad -  \phi^T(L_B\lambda  - CD_B\hat C^{ T}\pi - CD_B(\rho^+ - \rho^-))  \label{eq:L(x,sigma)}
\end{align}
where $\Phi_{S}(\lambda_S,\nu_S):=\sum_{i\in S}\Phi_i(\lambda_i,\nu_i)$ and
\begin{align}\label{eq:Phi_i}
\Phi_i(\lambda_i,\nu_i):=& c_i(d_i(\lambda_i+\nu_i)) - (\lambda_i+\nu_i)d_i(\lambda_i+\nu_i)\nonumber\\
&-D_i\frac{\nu_i^2}{2} + (\lambda_i+\nu_i)P_i^{in}.
\end{align}
The strict convexity of $L(d,\omega,x,\sigma)$ on $(d,\omega)$ and the fact that $d$ and $\omega$ only appear in \eqref{eq:olfc-b} and \eqref{eq:olfc-c} gives rise to the following lemma whose proof is also in the Appendix.
\begin{lemma}[Strict concavity of $L(x,\sigma)$ in $(\lambda,\nu)$]\label{lem:strict-concavity}
The function $\Phi_i(\lambda_i,\nu_i)$  in \eqref{eq:Phi_i} is strictly concave. As a result, $L(x,\sigma)$ is strictly concave in $(\lambda,\nu)$.
\end{lemma}

Next we reduce the Lagrangian $L(x,\sigma)$ by maximizing it for any $\nu_i$ with $i\in \L$. We let $y =\allowbreak (\lambda,\allowbreak\nu_\G,\pi,\rho^+,\rho^-)$ and consider the Lagrangian
\begin{equation}
L(x,y) = \underset{\nu_i:i\in\L}{\text{\textnormal{maximize}}}\; L(x,\sigma).\label{eq:reduced-lagrangian}
\end{equation}
Since $L(x,\sigma)$ is strictly concave in $\nu$ by Lemma \ref{lem:strict-concavity}, the minimizer of \eqref{eq:reduced-lagrangian} is unique. Moreover, this also implies that $L(x,y)$ is strictly concave in $(\lambda,\nu_\G)$.

Finally, the optimal load controllers can be then obtained by considering the primal-dual gradient law of $L(x,y)$ which is given by
\begin{align}\label{eq:red-primal-dual}
\dot y  =Y\left[\frac{\partial}{\partial y } L(x,y )^T\right]^+_\rho\;\text{ and }\;
\dot x  =- X\frac{\partial}{\partial x} L(x,y )^T
\end{align}
where $Y=\diag((\zeta^\nu_i)_{i\in\G},(\zeta^\lambda_i)_{i\in\N},(\zeta^\pi_k)_{k\in\K}, (\zeta^{\rho^+}_e)_{e\in\E},\allowbreak(\zeta^{\rho^-}_e)_{e\in\E})$ and $X=\diag((\chi^P_e)_{e\in\E},(\chi^\phi_i)_{i\in\N})$, and the projection $[\cdot]^+_\rho$ is defined as in \eqref{eq:projection} for $u=y$ and $u_S=\rho$.
This projection ensures that the vector $\rho(t)$ remains within the positive orthant, that is $\rho^+(t)\geq0$ and $\rho^-(t)\geq0$ $\forall t$.

\color{black}

The following theorem shows that this procedure indeed embeds the network dynamics as part of the primal-dual law \eqref{eq:red-primal-dual} while providing a distributed scheme to solve OLC.

\begin{theorem}[Optimal Load-side Control]\label{lem:primal-dual-gradient-law}
By setting $\zeta^\nu_i=M_i^{-1}$, $\chi^P_{ij}=B_{ij}$ and $\nu_i=\omega_i$, the primal-dual gradient law \eqref{eq:red-primal-dual} is equivalent to the power network dynamics \eqref{eq:1} together with 
\begin{subequations}\label{eq:load-control}
\begin{align}
\dot\lambda_i &=\zeta_i^\lambda\Big(P^{in}_i -d_i -\textstyle\sum_{j\in\N_i}B_{ij}(\phi_i-\phi_j)\Big)		\label{eq:load-control-a}\\
\dot\pi_k &= \zeta^\pi_k\Big(\textstyle\sum_{ij\in\B_k}\hat C_{k,ij}B_{ij}(\phi_i-\phi_j)-\hat P_{k}\Big)		\label{eq:load-control-b}\\
\dot\rho_{ij}^+ &= \zeta_{ij}^{\rho^+}\big[B_{ij}(\phi_i-\phi_j) - \overline{P}_{ij}\big]^+_{\rho_{ij}^+}							\label{eq:load-control-c}\\
\dot\rho_{ij}^- &= \zeta^{\rho_{ij}^-}\big[\underline P_{ij} -B_{ij}(\phi_i-\phi_j) \big]^+_{\rho_{ij}^-}							\label{eq:load-control-d}\\
\dot \phi_i&= \chi_i^\phi\Big(\textstyle\sum_{j\in\N_i} B_{ij}(\lambda_i\!-\!\lambda_j)\!-\! \textstyle \sum_{k\in\K,e\in\B_k}C_{i,e}B_e\hat C_{k,e}\pi_k 		\nonumber\\
				& \qquad\quad- \textstyle\sum_{e\in\E}C_{i,e}B_e(\rho_e^{+}-\rho_e^{-}) \Big) 																					\label{eq:load-control-e}\\
d_i &= {c_i'}^{-1}(\lambda_i+\omega_i) 																									\label{eq:load-control-f}
\end{align}
\end{subequations}
where \eqref{eq:load-control-a}, \eqref{eq:load-control-e} and \eqref{eq:load-control-f} are for $i\in\N$, \eqref{eq:load-control-b} is for $k\in\K$, and \eqref{eq:load-control-c} and \eqref{eq:load-control-d} are for $ij\in \E$.
\end{theorem}

\begin{IEEEproof}
By Lemma \ref{lem:strict-concavity} and \eqref{eq:L(x,sigma)}, $L(x,\sigma)$ is strictly concave in $(\lambda,\nu)$. Therefore, it follows that there exists a unique 
\begin{equation}\label{eq:nu_L^*}
\nu^*_\L(x,y) = \argmax_{\nu_\L} L(x,\sigma).
\end{equation} 
Moreover, by stationarity,  $\nu_\L^*(x,y)$ must satisfy
\begin{subequations}\label{eq:implicit-condition}
\begin{align}
&\frac{\partial L}{\partial \nu_\L} (x,y,\nu^*_\L(x,y))^T=\frac{\partial \Phi_\L}{\partial \nu_\L} (\lambda_\L,\nu^*_\L(x,y))^T - C_\L P\\
&=P^{in}_\L \!-\!D_\L\nu^*_\L(x,y) \!-\! d_\L(\lambda_\L\!+\!\nu^*_\L(x,y)) \!-\! C_\L P= 0
\end{align}
\end{subequations}
which is equivalent to \eqref{eq:1b}, i.e. $\nu^*_\L(x,y)$ implicitly satisfies \eqref{eq:1b}. %

We now iteratively apply the envelope theorem~\cite{mas1995microeconomic} on $L(x,y)$
\eqref{eq:reduced-lagrangian} to compute $\ddx{x}L(x,y)$ and $\ddx{y}L(x,y)$.
For example, to compute  $\ddx{x}L(x,y)$ we use
\begin{subequations}\label{eq:red-pLpx}
\begin{align}
&\ddx{x}L(x,y) = \left.\ddx{x}L(x,\sigma)\right\rvert_{\nu_\L=\nu^*_\L(x,y)}\\
&\!=\!\left. \left(\left. \ddx{x}L(d,\omega,x,\sigma)\right\rvert_{(d,\omega)=(c'^{-1}(\lambda+\nu),\nu)} \right) \right\rvert_{\nu_\L=\nu^*_\L(x,y)},
\end{align}
\end{subequations}
where $c'^{-1}(\lambda+\nu):=(c_i'^{-1}(\lambda_i+\nu_i))_{i\in\N}$, which leads to 
\begin{subequations}\label{eq:pLp-x}
\begin{align}
\frac{\partial}{\partial P}L(x,y) ^T&= -(C_\G^T\nu_\G + C_\L^T\nu^*_\L(x,y))								\label{eq:pLp-P}\\
\frac{\partial}{\partial \phi}L(x,y)^T&= -(L_B\lambda -CB(\hat C^{ T}\pi+ \rho^+-\rho^-))				\label{eq:pLp-phi}
\end{align}
\end{subequations}
An analogous computation for $\ddx{y}L(x,y)$ gives
\begin{subequations}\label{eq:pLp-y}
\begin{align}
\frac{\partial}{\partial \nu_\G}L(x,y)^T&=P^{in} -(d_\G(\lambda_\G+\nu_\G) +D\nu_\G) -CP 				\label{eq:pLp-nu_G}\\
\frac{\partial}{\partial \lambda}L(x,y)^T&=P^{in} -\left.d(\lambda+\nu)\right\rvert_{\nu_\L=\nu^*_\L(x,y)} -L_B\phi														\label{eq:pLp-lambda}\\
\frac{\partial}{\partial \pi}L(x,y)^T &= \hat C BC^T\phi-\hat P			\label{eq:pLp-pi}\\
\frac{\partial}{\partial \rho^+}L(x,y) ^T&= BC^T\phi - \overline P						\label{eq:pLp-rho+}\\
\frac{\partial}{\partial \rho^-}L(x,y) ^T&= \underline P - D_BC^T\phi 						\label{eq:pLp-rho-}
\end{align}
\end{subequations}
where, for a set $S$, we $d_S(\lambda_S+\nu_S):=(d_i(\lambda_i+\nu_i))_{i\in S}$ and $d(\lambda+\nu)=d_\N(\lambda_\N+\nu_\N)$. 

Using \eqref{eq:pLp-x}-\eqref{eq:pLp-y} and identifying $\omega_i$ in \eqref{eq:1} and \eqref{eq:load-control} with $\nu_i$ in \eqref{eq:red-primal-dual}, it is easy to see that \eqref{eq:red-primal-dual} only differs from \eqref{eq:1a},\eqref{eq:1c} and \eqref{eq:load-control} on the locations where $\omega_\L$ must be substituted with  $\nu_\L^*(x,y)$. 
However,  since there is a unique $\omega_\L$ that satisfies \eqref{eq:1b} given the remaining state variables, and we showed that $\nu^*_\L(x,y)$ also satisfies \eqref{eq:1b}, then it follows that \eqref{eq:red-primal-dual} and \eqref{eq:1} with \eqref{eq:load-control} are equivalent representations of the same system.
\end{IEEEproof}

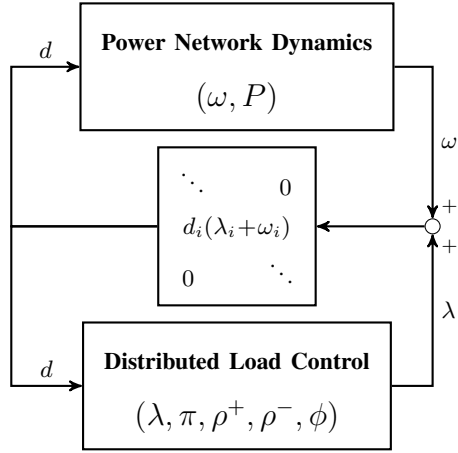
\begin{figure}
\centering
\begin{tikzpicture}[scale=1, auto,
 	decision/.style = { diamond, draw=blue, thick, fill=blue!20,
                        text width=5em, text badly centered,
                        inner sep=1pt, rounded corners },
    line/.style     = { draw, thick, ->, shorten >=2pt },
    myncbar/.style = {to path={
        let
            \p1=($(\tikztotarget)+(\myangle:\myarm)$)
        in
            -- ++(\myangle:\myarm) coordinate (tmp)
            -- ($(\tikztotarget)!(tmp)!(\p1)$)
            -- (\tikztotarget)\tikztonodes
    }},
    >=stealth']
    \small
    \matrix[ampersand replacement=\&, row sep=0.2cm, column sep=0.4cm] {
      \&
      \node[block] (pn) {\large $\begin{array}{c}
      								\textbf{\small Power Network Dynamics} \\
								      (\omega,P)
								     \end{array} $}; 
      \& \\

      \&
     \node[block] (L1) {$\begin{array}{ccc}
     	\ddots\!\!\!\!\!\!		& 	 				&		\!\!\!\!\!\!\!\! 0 \!\!\!\!		\\ 										& \!\!\!\!\!\! \small d_i(\lambda_i\!+\!\omega_i) \!\!\!\!\!\!	& 						\\
	\!\!\!\! 0 \!\!\!\!\!\!\!\!	&					& 		\!\!\!\!\!\!\ddots
        \end{array}$};\&
      \node [sum] (e1) {}; \\

	 \&
      \node[block] (lc) { \large $\begin{array}{c}
      									\textbf{\small Distributed Load Control}\\
										(\lambda,\pi,\rho^+,\rho^-,\phi)
				          			\end{array} $};
      \& \\
      
      \node[guide] (i1) {}; 
      \& \& \\
    };

    \draw [connector] (pn) -| node[pos=0.96] {\scriptsize$+$} node[near end] {$\omega$} (e1);
    \draw [connector] (e1) -- (L1);
    \draw [connector] (lc) -| node[pos=0.96, swap] {\scriptsize$+$} node [near end, swap] {$\lambda$} (e1);
	\draw [->,thick] (L1) to[myncbar,angle=180,arm=3]  node[auto] {$d$} (lc);
	\draw [->,thick] (L1) to[myncbar,angle=180,arm=3]  node[auto] {$d$} (pn);
  \end{tikzpicture}
  \caption{Control architecture derived from OLC}\label{fig:control-architecture}
\end{figure}

\begin{figure}[h!]\centering
\includegraphics[width=.425\columnwidth]{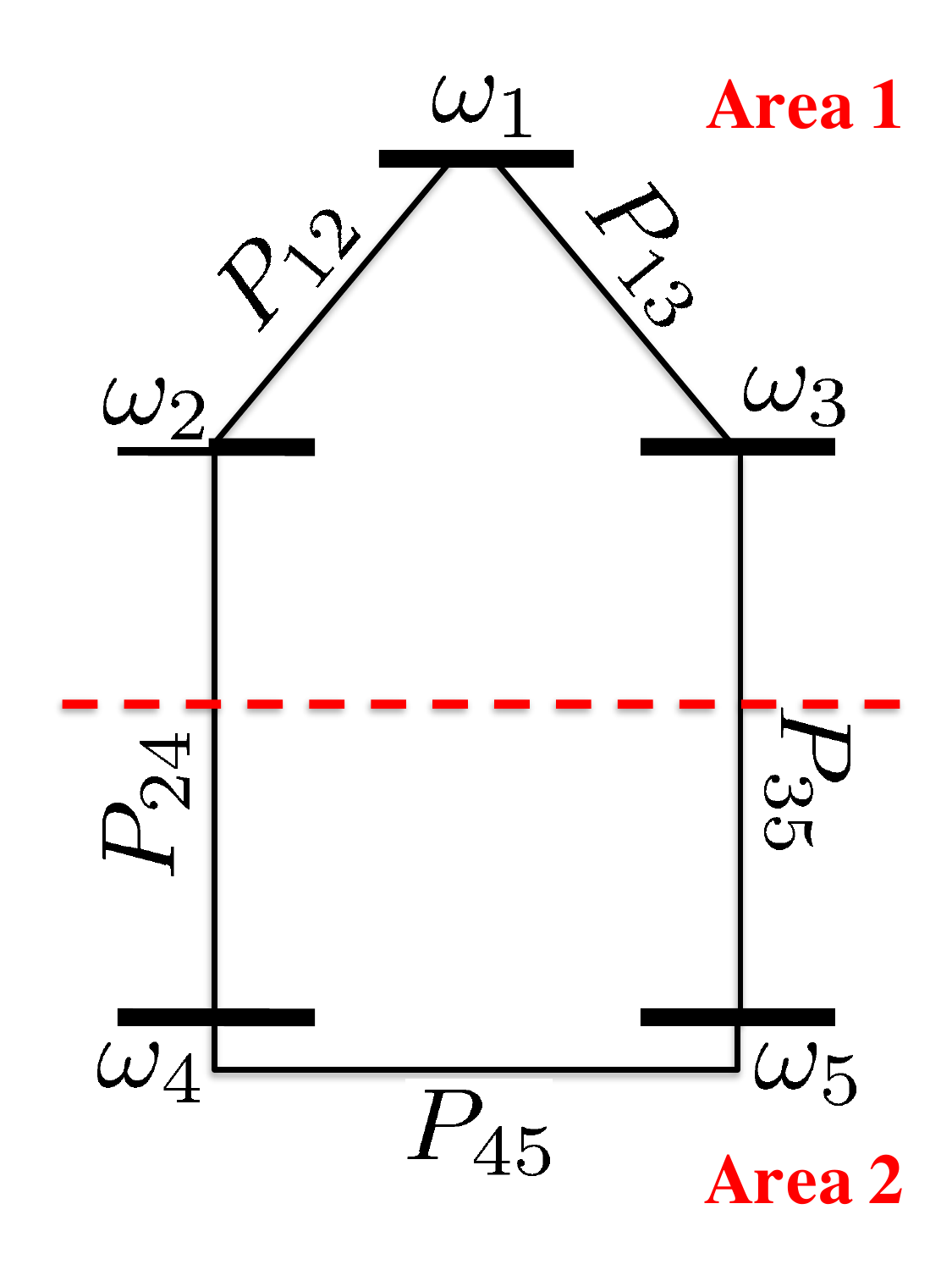}
\includegraphics[width=.425\columnwidth]{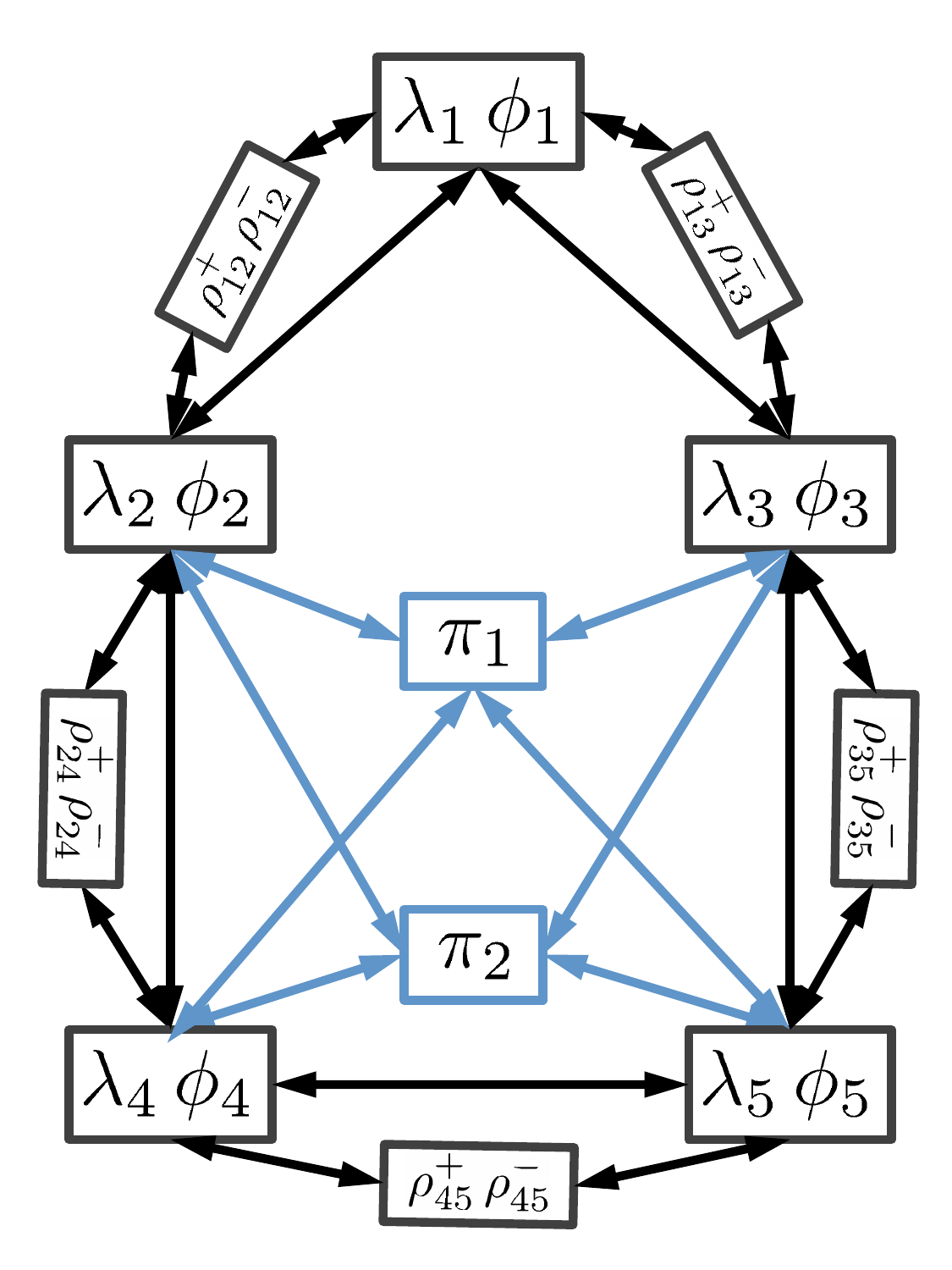}
\caption{Power network example (left) and the corresponding communication requirement to implement the distributed load control \eqref{eq:load-control}}\label{fig:diagram}
\end{figure}

Equations \eqref{eq:1} and \eqref{eq:load-control} show how the network dynamics can be complemented with dynamic load control such that the whole system amounts to a distributed primal-dual algorithm that tries to find a saddle point on $L(x,y)$. We will show in the next section that this system does achieve optimality as intended.

Figure \ref{fig:control-architecture} shows the unusual control architecture derived from our OLC problem. Unlike traditional observer-based controller design architecture~\cite{zhou1996robust}, our dynamic load control block does not try to estimate the state of the network. Instead, it drives the network towards the desired state using a {\it shared} static feedback loop, i.e. $d_i(\lambda_i+\omega_i)$. 

\added{Finally, Figure \ref{fig:diagram} shows the communication requirements to implement our distributed load control \eqref{eq:load-control}. The only state that cannot be computed distributedly is $\pi_k$ as it requires information from all boundary buses of area $k$ and their adjacent buses outside $k$. This limitation is overcome in Section \ref{ssec:inter-area-flow-constraints} by modifying \eqref{eq:welfare-d} in Problem~\ref{problem:olc}.}


\begin{remark}\label{rem:implementation-problem}
One of the limitations of \eqref{eq:load-control} is that in order to generate the Lagrange multipliers $\lambda_i$ one needs to estimate $P^{in}_i-d_i$ which is not easy to obtain from the measurements of $P^{in}_i-D_i\omega_i-d_i$ without knowing $D_i$. This problem will be addressed in Section \ref{sec:implementation-details} where we propose a modified control scheme that can achieve the same equilibrium without needing to know $D_i$ exactly.
\end{remark}
\begin{remark}
The procedure described in this section is independent of the constraints \eqref{eq:welfare-d}-\eqref{eq:welfare-e}. Therefore, such constraints can be generalized to arbitrary equality and inequality constraints on the line flows $BC^T\theta$. This property will be exploited in Section \ref{sec:framework-extensions} to further extend our framework.
\end{remark}


\section{Optimality and Convergence}
\label{sec:analysis}
In this section we will show that the system \eqref{eq:1} and \eqref{eq:load-control} can efficiently rebalance supply and demand, restore the nominal frequency, and preserve inter-area flow schedules and thermal limits. 

We will achieve this objective in two steps. Firstly, we will show that every equilibrium point of \eqref{eq:1} and \eqref{eq:load-control} is an optimal solution of \eqref{eq:olfc}, or equialently \eqref{eq:welfare}. This guarantees that a stationary point of the system efficiently balances supply and demand and achieves zero frequency deviation.

Secondly, we will show that every trajectory $(d(t),\allowbreak \omega(t),\allowbreak \phi(t),P(t),\lambda(t),\pi(t),\rho^+(t),\rho^-(t))$ converges to an equilibrium point of \eqref{eq:1} and \eqref{eq:load-control}. Moreover, we will show that since $P(0)=BC^T\theta^0$ (as shown in Section \ref{ssec:network-model}), the line flows will converge to a point that satisfies~\eqref{eq:inter-area-constraint-2} and \eqref{eq:thermal-limit}. 
\begin{theorem}[Optimality]\label{thm:optimality}
A point $p^*:=(d^*,\omega^*,\phi^*,P^*,\lambda^*,\allowbreak\pi^*,\rho^{+*},\rho^{-*})$ is an equilibrium point of  \eqref{eq:1} and \eqref{eq:load-control} if and only if $(d^*,\omega^*,\theta^*)$ is an optimal solution of OLC and $(d^*,\omega^*,\phi^*,P^*,\lambda^*,\nu^*,\pi^*,\rho^{+*},\rho^{-*})$ is a primal-dual optimal solution to the VF-OLC problem, with 
\begin{equation}\label{eq:condition-thm}
\omega^*=\nu^*,\text{ }CP^*=L_B\theta^*\text{ and } C^T\theta^*=C^T\phi^*.
\end{equation}
\end{theorem}
\begin{IEEEproof}
The proof of this theorem is a direct application of lemmas \ref{lem:OLC-optimality} and \ref{lem:equivalence}. 
Let $p^*$ be an equilibrium of \eqref{eq:1} and \eqref{eq:load-control}.
Then, by definition of the projection $[\cdot]_{\rho}^+$ and \eqref{eq:load-control-c}-\eqref{eq:load-control-d}, $\rho^{+*}\geq 0$ and $\rho^{-*}\geq 0$ and thus dual feasible.

Similarly, since $\dot\omega_i=0$, $\dot\lambda_i=0$, $\dot\pi_{k}=0$, $\dot\rho_{ij}^+=0$ and $\dot\rho_{ij}^-=0$, then \eqref{eq:1a}-\eqref{eq:1b} and \eqref{eq:load-control-a}-\eqref{eq:load-control-d} are equivalent to primal feasibility, i.e. $(d^*,\omega^*,\phi^*,P^*)$ is a feasible point of \eqref{eq:olfc}.
Finally, by definition of \eqref{eq:1} and \eqref{eq:load-control}, conditions \eqref{eq:16}, \eqref{eq:16b} and \eqref{eq:complementary-slackness} are always satisfied by any equilibrium point. Thus we are under the conditions of Lemma \ref{lem:OLC-optimality} and therefore $(d^*,\omega^*,\phi^*,P^*,\lambda^*,\nu^*,\allowbreak\pi^*,\rho^{+*},\rho^{-*})$ is primal-dual optimal of VF-OLC satisfying \eqref{eq:condition-thm}. Lemma \ref{lem:equivalence} shows the remaining statement of the theorem.
\end{IEEEproof}


The rest of this section is devoted to showing that in fact for every initial condition $(\omega(0),\phi(0),P(0),\lambda(0),\pi(0),\allowbreak\rho^+(0),\rho^-(0))$, the system \eqref{eq:1} and \eqref{eq:load-control} converges to one such optimal solution. Furthermore, we will show that $P(t)$ converges to a $P^*$ that satisfies \eqref{eq:inter-area-constraint-2} and \eqref{eq:thermal-limit}.

Since we showed in Theorem \ref{lem:primal-dual-gradient-law}  that \eqref{eq:1} and \eqref{eq:load-control} is equivalent to \eqref{eq:red-primal-dual}, we will provide our convergence result for \eqref{eq:red-primal-dual}. 
Our global convergence proof builds on recent results of \cite{Feijer:2010ve} on global convergence of primal-dual algorithms for network flow control. Our proof extends \cite{Feijer:2010ve} in the following aspects. Firstly, the Lagrangian $L(x,y)$ is not strictly concave in all of its variables. Secondly, the projection \eqref{eq:projection} introduces discontinuities in the vector field that prevents the use of the standard LaSalle's Invariance Principle~\cite{khalil2002nonlinear}.

We solve the latter issue using an invariance principle for Caratheodory systems~\cite{Bacciotti:2006uj}. We refer the reader to~\cite{cherukuri2015Convergence} for a detailed treatment that formalizes its application for primal-dual systems.
The former issue is solved in Theorem \ref{thm:global-convergence} which makes use of the following additional lemma whose proof can be found in the Appendix.

\begin{lemma}[Differentiability of $\nu^*_\L(x,y)$]\label{lem:differentiability-nu}
Given any $(x,y)$, the maximizer of \eqref{eq:reduced-lagrangian}, $\nu^*_\L(x,y)$, is continuously differentiable provided
$c_i(\cdot)$ is strongly convex. Furthermore, the derivative is given by
\begin{align}
\frac{\partial}{\partial x}\nu_\L^*(x,y)\! &= \!\!
\begin{blockarray}{c@{\hspace{2.5pt}}cc@{\hspace{2.5pt}}cl}
  & \matindex{\phi} & \matindex{P} & &\\
\begin{block}{[c@{\hspace{2.5pt}}c|c@{\hspace{2.5pt}}c]l}
 & \! 0 \!&\! -(D_\L+ d'_\L)^{-1}C_\L  \! & &\! \matindex{\nu_\L}\\
\end{block}
\end{blockarray}\label{eq:dnuLdx}
\\
\frac{\partial}{\partial y}\nu^*_\L(x,y) &\!=\!\!
\begin{blockarray}{c@{\hspace{2.5pt}}ccccc@{\hspace{2.5pt}}cl}
  &\!\!\matindex{\lambda_\L}\!\! &  \!\!\matindex{\lambda_\G}\!\! & \!\!\matindex{\nu_\G}\!\!  & \!\!\matindex{\pi}\!\!  & \!\!\matindex{\rho}\!\! & &\\
\begin{block}{[c@{\hspace{2.5pt}}c|c|c|c|c@{\hspace{2.5pt}}c]l}
 & \!\!-(D_\L\!+\! d'_\L)^{-1} d'_\L\!\! & \!\!0\!\! & \!\!0\!\! & \!\!0\!\! & \!\!0\!\!  & &\! \matindex{\nu_\L}\\
\end{block}
\end{blockarray}\label{eq:dnuLdy}
\end{align}
where $D_S:=\diag(D_i)_{i\in S}$ and 
$$
d'_S = 
\begin{cases}
\diag(d'_i(\lambda_i+\nu_i^*(x,y)))_{i\in S} & \text{ if } S\subseteq \L\\
\diag(d'_i(\lambda_i+\nu_i))_{i\in S} & \text{ if } S\subseteq \G\\
\end{cases}
$$
\end{lemma}

We now present our main convergence result. Let $E$ be the set of equilibrium points of \eqref{eq:red-primal-dual}
\[
E:=\left\{(x,y): \frac{\partial L}{\partial x}(x,y)=0,\;\left[\frac{\partial L}{\partial y}(x,y)\right]^+_{\rho}=0\right\},
\]
which by theorems \ref{lem:primal-dual-gradient-law} and \ref{thm:optimality} characterizes the set of optimal solutions of the OLC problem.
\begin{theorem}[Global Convergence]\label{thm:global-convergence}
The set $E$ of equilibrium points of the primal-dual algorithm \eqref{eq:red-primal-dual} is
 globally asymptotically stable. Furthermore, each individual trajectory converges to a unique point within $E$ that is optimal with respect to the OLC problem.
\end{theorem}
\begin{IEEEproof}
Following \cite{Feijer:2010ve} we consider the candidate Lyapunov function
\begin{align}\label{eq:U}
U(x,y)&=\! \frac{1}{2}(x\!-\!x^*)^TX^{-1}(x\!-\!x^*) 
\!+\! \frac{1}{2}(y-y^*)^TY^{-1}(y\!-\!y^*)
\end{align}
where $(x^*,y^*)$ is {\it any} equilibrium point of \eqref{eq:red-primal-dual}.

\color{black}
We divide the proof of this theorem in four steps:
\begin{enumerate}[leftmargin=*,label=\bfseries Step \arabic*:, wide=\parindent]
\item We first use the invariance principle for Caratheodory systems~\cite{Bacciotti:2006uj} to show that $(x(t),y(t))$ converges to the largest invariance set that satisfies $\dot U(x,y) \equiv 0$  between transitions of the projection $[\cdot]^+_\rho$, i.e. 
\begin{align}
(x(t),y(t))\rightarrow M\subseteq\{(x,y)\!:\!  \dot U(x(t),y(t))\!\equiv\! 0, t\in \mathds R^+\backslash\{t_k\}\}\label{eq:M}
\end{align}
where $\{t_k,\,k\in\mathds N\}$ are the time instants when the projection changes between on and off.
\item We show that any invariant trajectory  $(x(t),y(t))\in M$ must have $\lambda(t)\equiv\hat\lambda$ and $\nu(t)\equiv \hat \nu$ for some constant vectors $\hat\lambda$ and $\hat \nu$.
\item We show that whenever $\lambda(t)\equiv\hat\lambda$ and $\nu(t)\equiv \hat \nu$, then whole trajectory $(x(t),y(t))$ must an equilibrium point, i.e. $M\subseteq E$.
\item Finally, we show that even though the invariance principle guarantees only convergence to the set $E$. The convergence is always to some point within $E$, i.e. $(x(t),y(t))\rightarrow(x^*,y^*)\in E$.
\end{enumerate}
\vspace{1.5ex}
\color{black}
{\bf Proof of step 1:}
Differentiating $U$ over time gives
\begin{align}
&\dot U(x,y) \!=\!-\ddx{x}L(x,y)(x\!-\!x^*) \!+\! \left[\ddx{y}L(x,y)\right]^+_\rho\!\!\!(y\!-\!y^*)\label{eq:dotU-0}\\
&\leq \frac{\partial L}{\partial x}(x,y) (x^*-x) + \frac{\partial L}{\partial y}(x,y)(y-y^*)\label{eq:dotU-1}\\
&\leq L(x^*,y) - L(x,y) + L(x,y) - L(x,y^*)\label{eq:dotU-2}\\
&= L(x^*,y)  - L(x,y^*) \nonumber\\
&= \underset{\leq 0}{\underbrace{L(x^*,y)  - L(x^*,y^*)}}   + \underset{\leq0}{\underbrace{L(x^*,y^*) - L(x,y^*)}}\label{eq:dotU-3}
\end{align}
where \eqref{eq:dotU-0} follows from \eqref{eq:red-primal-dual} and \eqref{eq:dotU-1} from \eqref{eq:projection-property}. 
Step \eqref{eq:dotU-2} follows from convexity (resp. concavity) of $L(x,y)$ in $x$ (resp. $y$).  
Finally, equation \eqref{eq:dotU-3} follows from the saddle property of the equilibrium point $(x^*,y^*)$.

Therefore, since $U(x,y)$ is radially unbounded, the trajectories are bounded, and it follows from the invariance principle for Caratheodory systems~\cite{Bacciotti:2006uj} that $(x(t),y(t))\rightarrow M$, i.e.  \eqref{eq:M} holds. 
The steps 2 and 3 below basically characterize $M$.


\vspace{1.5ex}
{\bf Proof of step 2:}
Notice that in order to have $\dot U \equiv 0$, both terms in \eqref{eq:dotU-3} must be zero. 
In particular, we must have 
\[
L(x(t),y^*)\equiv L(x^*,y^*).
\]
Now, differentiating with respect to time gives
\begin{align*}
0&\equiv\dot{\overline{L(x(t),y^*)}}\equiv \ddx{x}L(x(t),y^*)\dot x
\equiv -||\ddx{x}L(x(t),y^*)||^2,
\end{align*}
which implies that $\ddx{x}L(x(t),y^*)\equiv 0.$


Therefore, the fact that $\nu_\G^*=0$, $\ddx{P}L(x(t),y^*)\equiv 0$, and \eqref{eq:pLp-P} holds, implies that $x(t)$ must satisfy $C_\L^T\nu_\L^*(x(t),y^*)\equiv 0$, which implies that either $\nu_\L^*(x(t),y^*)\equiv 0$ (when $C_\L$ is full row rank) or
 $\nu_\L^*(x(t),y^*)\equiv 1\alpha(t)$ (when $\L=\N$) where $\alpha(t)$ is a time-varying scalar.

We now show that when $\L=\N$ we get $\nu_\L^*(x(t),y^*)\equiv \hat \nu_\L$ for some constant vector $\hat\nu_\L$. Differentiating $\nu_\L^*(x(t),y^*)\equiv 1\alpha(t)$ with respect to time and using \eqref{eq:dnuLdx} we obtain
\[
(D_\L+d'_\L)^{-1}C_\L \dot P((t) \equiv \1 \dot \alpha(t)
\]
which after left multiplying by $\1^T(D_\L+d'_\L)$ gives
\[
\1^T(D_\L+d'_\L)\1 \dot \alpha(t)\equiv 0\implies \dot \alpha(t)\equiv 0.
\]
Thus,  in either case we obtain
\begin{equation}\label{eq:th9-1}
\nu_\L^*(x(t),y^*)=\nu_\L^*(C_\L P(t),\lambda_\L^*)\equiv \hat \nu_\L
\end{equation}
for some constant vector $\hat \nu_\L$, which implies that $C_\L P(t)\equiv C_\L\hat P$ for some constant vector $\hat P$.

Therefore, it follows that $\nu^*_\L(x(t),y(t))$ must satisfy
\begin{align}\label{eq:th9-2}
\nu^*_\L(x(t),y(t))\equiv\nu^*_\L(\hat x,y(t))
\end{align}
for some constant vector $\hat x$.

Now, using \eqref{eq:implicit-condition} with \eqref{eq:th9-2} we get
\begin{align}\label{eq:th9-22}
&P^m _\L- D_\L\nu^*_\L(\hat x,y(t)) -d_\L(\lambda_\L(t)+\nu^*_\L(\hat x,y(t))) -C_\L \hat P \nonumber\\
&\equiv0.
\end{align}

A similar argument using the fact that $L(x^*,y)\equiv L(x^*,y^*)$ gives
\begin{equation}\label{eq:th9-proj}
\ddx{y}L(x^*,y)\left[\ddx{y}L(x^*,y)^T\right]_\rho^+\equiv 0.
\end{equation}
Since the projection $[\cdot]^+_\rho$ only acts on the $\rho$ positions, \eqref{eq:th9-proj} implies $\ddx{\nu_\G}L(x^*,y)\equiv 0$, $\ddx{\lambda}L(x^*,y)\equiv 0$ and $\ddx{\pi}L(x^*,y)\equiv 0$.

Now $\ddx{\nu_\G}L(x^*,y)\equiv 0$ together with  equation \eqref{eq:pLp-nu_G} implies that
\begin{align}
P^m_\G - D_\G\nu_\G (t) -d_\G(\lambda_\G(t)+\nu_\G(t)) -C_\G P^*\equiv 0,\label{eq:th9-3}
\end{align}
and $\ddx{\lambda}L(x^*,y)\equiv 0$ with \eqref{eq:pLp-lambda} implies
\begin{align}
P^m_\G  &-d_\G(\lambda_\G(t)+\nu_\G(t)) -C_\G P^*\equiv 0\label{eq:th9-4}\\
P^m_\L  &-d_\L(\lambda_\L(t)+\nu^*_\L(x^*,y(t))) -C_\L P^*\equiv 0\label{eq:th9-5}
\end{align}

Using \eqref{eq:th9-3} and \eqref{eq:th9-4} together with the fact that $d_i(\cdot)$ is strictly increasing, we get
$\nu_\G(t)\equiv\hat \nu_\G$ and $\lambda_\G(t)\equiv\hat\lambda_\G$, for constant vectors $\hat \nu_\G$ and $\hat \lambda_\G$. Moreover, since $P^*$ is primal optimal, Lemma \ref{lem:primal-dual-gradient-law} and Theorem \ref{thm:optimality} imply that $\nu_\G(t)\equiv 0$ and $\lambda_\G(t)\equiv \lambda_\G^*$. Finally, now using \eqref{eq:th9-22} together with \eqref{eq:th9-5}, the same argumentation gives $\nu^*_\L(x(t),y(t)) \equiv \hat\nu_\L$ and $\lambda_\L(t)\equiv\hat \lambda_\L$ for constant vectors $\hat\nu_\L$ and $\hat \lambda_\L$.
This finishes step 2, i.e. $\lambda(t)\equiv \hat \lambda$ and $\nu(t)\equiv \hat\nu$.

\vspace{1.5ex}
{\bf Proof of step 3:}
Now, since $\dot\lambda\equiv 0$, it follows from \eqref{eq:load-control-a} that $C^T\phi(t)\equiv C^T\hat \phi$ for some constant vector $\hat \phi$ or equivalently $\phi(t)\equiv \hat \phi + \beta(t)\1$. Differentiating in time $\1^T(\chi^\phi)^{-1}\phi(t)$ gives $0\equiv\1^T(\chi^\phi)^{-1}\dot \phi\equiv(\sum_{i\in\N}{1}/{\chi^\phi_i})\dot\beta$ which implies that $\beta(t)\equiv \hat\beta$ for constant scalar $\hat\beta$.

Suppose now that either $\dot P\neq 0$ or $\dot \pi\neq 0$. Since $C^T\phi(t)\equiv C^T\hat \phi$ and $\nu(t)\equiv\hat\nu$,   $\dot P$ and $\dot \pi$ are constant.
Thus, since the trajectories are bounded, we must have $\dot P\equiv0$ and $\dot \pi\equiv0$; otherwise $U(x,y)$ will grow unbounded (contradiction). 

It remains to show that $\dot \rho\equiv0$, i.e. $\dot \rho^+\equiv\dot\rho^-\equiv0$. 
Since $\phi(t)\equiv \hat \phi$, then the argument inside \eqref{eq:load-control-c} and \eqref{eq:load-control-d} is constant. 

Now consider any $\rho_e^+$, $e=ij\in\E$. Then we have three cases:
(i) $B_{e}(\hat \phi_i - \hat \phi_j) - \bar P_e >0$, (ii) $B_e(\hat \phi_i - \hat \phi_j) - \bar P_e <0$ and (iii) $B_e(\hat \phi_i - \hat \phi_j) - \bar P_e =0$.
Case (i) implies $\rho^+_e(t)\rightarrow +\infty$ which cannot happen since the trajectories are bounded. Case (ii) implies that $\rho^+_e(t)\equiv 0$ which implies that $\dot\rho^+_e\equiv0$, and case (iii) also implies  $\dot\rho^+_e\equiv0$.
An analogous argument gives $\dot\rho^-\equiv0$.
Thus, we have shown that $M\subseteq E$.

\vspace{1.5ex}
{\bf Proof of step 4:}
We now use structure of $U(x,y)$ to achieve convergence to a single equilibrium. Since $(x(t),y(t))\rightarrow M$ and $(x(t),y(t))$ are bounded, then there exists an infinite  sequence of time values $\{t_k\}$ such that $(x(t_k),y(t_k))\rightarrow (\hat x^*,\hat y^*)\in M$. Thus, using this specific equilibrium $(\hat x^*,\hat y^*)$ in the definition of $U(x,y)$, it follows that  $U(x(t_k),y(t_k))\rightarrow 0$, which by continuity of $U(x,y)$ implies that $(x(t),y(t))\rightarrow(\hat x^*,\hat y^*)$.

Thus, it follows that $(x(t),y(t))$ converges to only one optimal solution within $M\subseteq E$.
\end{IEEEproof}

Finally, the following theorem shows that the system is able to restore the inter-area flows~\eqref{eq:inter-area-constraint-2} and maintain the line flows within the thermal limits~\eqref{eq:thermal-limit}.
\begin{theorem}[Inter-area Constraints and Thermal Limits]\label{thm:inter-area-flows}
Given any primal-dual optimal solution $(x^*,\sigma^*)\in E$, the optimal line flow vector $P^*$ satisfies \eqref{eq:inter-area-constraint-2}. Furthermore, if 
$P(0) = BC^T\theta^0$ for some $\theta^0\in\R^{|\N|}$, then $P_{ij}^*=B_{ij}(\phi_i^*-\phi_j^*)$ and therefore \eqref{eq:thermal-limit} holds.
\end{theorem}
\begin{IEEEproof}
By optimality, $P^*$ and $\phi^*$ must satisfy
\begin{align}\label{eq:thm9-1}
P^m-d^*&= CP^* = L_B\phi^*=CBC^T\phi^*
\end{align}

Therefore using primal feasibility, \eqref{eq:Cflat} and \eqref{eq:thm9-1} we have
\begin{align*}
\hat P \!=\! \hat C BC^T\phi^*\!=\!E_\K CBC^T\phi^*\!=\!E_\K CP^*\!=\!\hat C P^*
\end{align*}
which is exactly \eqref{eq:inter-area-constraint-2}.

Finally, to show that $P_{ij}^*= B_{ij}(\phi_i^*-\phi_j^*)$ we will use 
\eqref{eq:1c}. Integrating \eqref{eq:1c}
over time gives
$$ P(t) - P(0) = \int_0^tBC^T \nu(s)ds.$$
Therefore, since $P(t)\rightarrow P^*$,  we have
 $P^* = P(0) + BC^T\theta^*$ 
where $\theta^*$ is any finite vector satisfying $C^T\theta^*=\int_0^\infty C^T\nu(s)ds$.

Again by primal feasibility $CBC^T\phi^*=L_B\phi^* = CP^* = C(P(0) + BC^T\theta^*)= CBC^T(\theta^0+\theta^*).$
Thus, we must have $\phi^*=(\theta^0+\theta^*)+\alpha \1$ and it follows then that $P^*=BC^T(\theta^0+\theta^*)=BC^T(\phi^*-\alpha\1)=BC^T\phi^*$. Therefore, since by primal feasibility $\underline P\leq BC^T\phi^*\leq \bar P$, then $\underline P\leq P^*\leq \bar P$.
\end{IEEEproof}

%
%


\section{Convergence under uncertainty}\label{sec:implementation-details}

In this section we discuss an important aspect of the implementation of the control law \eqref{eq:load-control}. We provide a modified control law that solves the problem raised in Remark \ref{rem:implementation-problem}, i.e. that does not require knowledge of $D_i$. We show that the new control law still converges to the same equilibrium provided the estimation error of $D_i$ is small enough (c.f. \eqref{eq:condition-th16}). 


We propose an alternative mechanism to compute $\lambda_i$. Instead of \eqref{eq:load-control-a}, we consider the following control law:
\begin{subequations}\label{eq:load-control-2}
\begin{align}
\dot\lambda_i &\!=\!\zeta_i^\lambda\Big(M_i\dot \omega_i\!+\!a_i\omega_i \!+\!\sum_{e\in\E}\!C_{i,e}P_e \!-\!\sum_{j\in\N_i}\!B_{ij}(\phi_i\!-\!\phi_j)\Big)						\label{eq:load-control-2-a}
\end{align}
\end{subequations}
where  $M_i:=0$ for $i\in\L$ and $a_i\in\R$ is a positive controller parameter that can be arbitrarily chosen.
Notice that, while before $D_i$ was an unknown quantity, $M_i$ is usually known and $a_i$ is a design parameter.
\added{Furthermore, while equation \eqref{eq:load-control-2-a} requires the knowledge of $\dot\omega_i$, this is only needed on generator buses and can therefore be measured from the generator's shaft angular acceleration
using one of several existing mechanisms, see e.g.~\cite{Tilli:2001wu}.}


The parameter $a_i$ plays the role of $D_i$. In fact, whenever $a_i=D_i$ then one can use \eqref{eq:1a}-\eqref{eq:1b} to show that \eqref{eq:load-control-2-a}  is the same as \eqref{eq:load-control-a}.
More precisely, if we let $a_i=D_i+\delta a_i$, then using \eqref{eq:1a}-\eqref{eq:1b}, \eqref{eq:load-control-2-a} becomes
\begin{align}
\dot\lambda 
&\!=\!\zeta_i^\lambda\Big( P_i^{in} \!-\!d_i \!+\!\delta a_i\omega_i \!-\!\sum_{j\in\N_i}B_{ij}(\phi_i\!-\!\phi_j)\Big),\label{eq:49}
\end{align}
which is equal to \eqref{eq:load-control-a} when $\delta a_i=0$.
\added{A simple equilibrium analysis shows that $a_i$ does not affect the steady state behavior provided that $a_i\not=0$ for some $i\in\N$. Thus, we focus in this section on studying the stability of our modified control law.
}

Using \eqref{eq:49}, we can express the new system using
\begin{subequations}\label{eq:red-primal-dual-2}
\begin{align}
\dot x  &=- X\frac{\partial}{\partial x} L(x,y )^T\\
\dot y  &=Y\left[\frac{\partial}{\partial y } L(x,y )^T+ g(x,y)\right]^+_\rho
\end{align}
\end{subequations}\vspace{-.15in}
\begin{flalign}\label{eq:th8-2}
&\!\!\!\!\!\!\!\!\!\text{where}\quad g(x,y)\!:=\!\!
\begin{blockarray}{c@{\hspace{2pt}}c@{\hspace{2pt}}c@{\hspace{2pt}}l}
\matindex{\lambda_\L}  &\matindex{\lambda_\G} & \matindex{\,(\nu_\G,\pi,\rho)\,} &\\
\begin{block}{[c@{\hspace{2pt}}c@{\hspace{2pt}}c@{\hspace{2pt}}]l}
(\delta A_\L\nu^*_\L)^T & (\delta A_\G\nu_\G)^T & 0 &\;^T \\
\end{block}
\end{blockarray}\vspace{-.2in}
\end{flalign}
with $\nu_\L^*:=\nu_\L^*(x,y)$ and $\delta A_S:=\diag(\delta a_i)_{i\in S}$.

The system \eqref{eq:red-primal-dual-2} is no longer a primal-dual algorithm. The main result of this section shows, that 
under provided that $a_i$ does not depart significantly from $D_i$ (see \eqref{eq:condition-th16}), convergence to the optimal solution is preserved.
%

To show this result, we provide a novel convergence proof that makes use of the following lemmas whose proofs can be found  in the Appendix.
\begin{lemma}[Second order derivatives of $L(x,y)$]\label{lem:second-order-derivatives}
Whenever Lemma \ref{lem:differentiability-nu} holds, then we have
\begin{align}
&\dddx{x}L(x,y)=
\begin{blockarray}{c@{\hspace{5pt}}cc@{\hspace{5pt}}cl}
  &\matindex{\phi} & \matindex{P} & &\\
\begin{block}{[c@{\hspace{5pt}}cc@{\hspace{5pt}}c]l}
& 0 &                                              0  & & \matindex{\phi}\\
&  0 & C_\L^T(D_\L+d'_\L)^{-1}C_\L  & & \matindex{P}\\
\end{block}
\end{blockarray}\label{eq:lem10-a}
\text{ and } \\
&\dddx{y}L(x,y)\!=\!-\!
\begin{blockarray}{c@{\hspace{1.5pt}}c@{\hspace{1.5pt}}c@{\hspace{1.5pt}}c@{\hspace{.5pt}}c@{\hspace{.5pt}}c@{\hspace{2pt}}l}
  &\matindex{\lambda_\L} &\matindex{\lambda_\G} & \matindex{\nu_\G} & \!\!\matindex{(\pi,\rho)}\!\! & &\\
\begin{block}{[c@{\hspace{1.5pt}}c@{\hspace{1.5pt}}c@{\hspace{1.5pt}}c@{\hspace{1.5pt}}c@{\hspace{1.5pt}}c@{\hspace{1.5pt}}]@{\hspace{2pt}}l}
& D_\L(D_\L+d'_\L)^{-1}d'_\L & 0 & 0 & 0 & & \,\matindex{\lambda_\L}\\
& 0 & d'_\G & d'_\G & 0 &  & \,\matindex{\lambda_\G}\\
& 0 & d'_\G &  (D_\G+d'_\G)  & 0 & & \,\matindex{\nu_\G}\\
& 0 & 0 & 0 & 0 & & \,\matindex{(\pi,\rho)} \\
\end{block}
\end{blockarray}\label{eq:lem10-b}
\end{align}
with $\dddx{x}L(x,y)\succeq 0$ and $\dddx{y}L(x,y)\preceq 0$.
\end{lemma}

\begin{lemma}[Partial derivatives of $g(x,y)$]\label{lem:partial-derivatives}
Whenever Lemma \ref{lem:differentiability-nu} holds, then
\begin{align*}
&\ddx{x}g(x,y)=
\begin{blockarray}{c@{\hspace{5pt}}cc@{\hspace{5pt}}cl}
  & \matindex{\phi} &\matindex{P} & &\\
\begin{block}{[c@{\hspace{5pt}}cc@{\hspace{5pt}}c]l}
& 0 & -\delta A_\L(D_\L+d'_\L)^{-1}C_\L & & \matindex{\lambda_\L}\\
& 0 & 0 & & (\matindex{\nu_\G,\lambda_\G},\pi,\rho)\\
\end{block}
\end{blockarray}
\\
&\ddx{y}g(x,y)=
\begin{blockarray}{c@{\hspace{.5pt}}c@{\hspace{.5pt}}c@{\hspace{.5pt}}c@{\hspace{.5pt}}c@{\hspace{.5pt}}cl}
  &\matindex{\lambda_\L} &\matindex{\lambda_\G} & \matindex{\nu_\G} & \matindex{(\pi,\rho)} & &\\
\begin{block}{[c@{\hspace{3.5pt}}c@{\hspace{3.5pt}}c@{\hspace{3.5pt}}c@{\hspace{3.5pt}}c@{\hspace{3.5pt}}c]l}
& -\delta A_\L(D_\L+d'_\L)^{-1}d'_\L & 0 & 0 & 0 & & \matindex{\lambda_\L}\\
& 0 & 0 & \delta A_\G & 0 & & \matindex{\lambda_\G}\\
& 0 & 0 & 0 & 0 & & \matindex{(\nu_\G,\pi,\rho)} \\
\end{block}
\end{blockarray}
\end{align*}
\end{lemma}

Unfortunately, the conditions of Theorem \ref{thm:global-convergence} will not suffice to guarantee convergence of the perturbed system. The main difficulty is that $d'_i(\lambda_i+\nu_i)>0$ can become arbitrarily close to zero. 
Therefore the sub-matrix of \eqref{eq:lem10-b} corresponding to the states $\lambda$ and $\nu_\G$ can become arbitrarily close to singular which makes the system not robust to perturbations of the form of \eqref{eq:th8-2}.
Thus, we require the following additional assumption.
\begin{assumption}[Lipschitz continuity of $c'_i$]\label{as:4}
The functions $c'_i(\cdot)$ are Lipschitz continuous with Lipschitz constant $L>0$.
\end{assumption}

\added{
Assumption \ref{as:4} implies that the domain of  $\D_i=\R$ in Assumption \ref{as:1}. However, if the systems is designed with enough capacity so that $d_i^*\in[\underline{d}_i +\varepsilon, \overline{d}_i-\varepsilon]$ $\forall i$, then
one can always modify a cost function $c_i(\cdot)$ that satisfies Assumptions \ref{as:1} and \ref{as:2} for finite domains $\D_i=[\underline{d}_i,\overline{d}_i]$ and define a new cost function $\tilde c_i(\cdot)$ that satisfies Assumption \ref{as:4} without modifying the optimal allocation $d_i^*$. More precisely, given $c_i(\cdot)$, define $\tilde c_i(\cdot)$ to be equal to $c_i(\cdot)$, inside $[\underline{d}_i +\varepsilon, \overline{d}_i-\varepsilon]$ and modify $\tilde c_i(\cdot)$ outside the subset so that Assumption \ref{as:4} holds.
It is easy to show then that the optimal solution will still be $d_i^*$ and therefore we still get $\underline{d}_i\leq d_i^*\leq \overline{d}_i$.
}

Using now Assumptions \ref{as:1} and \ref{as:4} we can show that $\alpha\leq c_i''\leq L$ which implies 
\begin{equation}\label{eq:d'-bounds}
\underline{d}':= {1}/{L} \quad \leq \quad d'_i={1}/{c_i''}\quad \leq \quad  \bar{d}':={1}/{\alpha}.
\end{equation}

\begin{theorem}[Global convergence of perturbed system]\label{thm:global-convergence-2}
Whenever assumptions \ref{as:1}, \ref{as:2} and \ref{as:4} hold. The system \eqref{eq:red-primal-dual-2} converges to a point in the optimal set $E$ for every initial condition whenever
\begin{equation}
\delta a_i\in2(\;\underline{d}' -\sqrt{{\underline{d}'}^2 +\underline{d}' D_{\min}},\; \underline{d}' +\sqrt{{\underline{d}'}^2 +\underline{d}' D_{\min}}\;).	\label{eq:condition-th16}
\end{equation}
where $D_{\min}:=\min_{i\in\N} D_i$.
\end{theorem}
\color{black}
\begin{IEEEproof}
We prove this theorem in three steps:
\begin{enumerate}[leftmargin=*,label=\bfseries Step \arabic*:, wide=\parindent]
\item We first show that under the dynamics \eqref{eq:red-primal-dual-2}, the time derivative of \eqref{eq:U} is upper-bounded by
\begin{equation}\label{eq:dotUz-integral}
\dot U(z) \leq \int_0^1 \!(z\!-\!z^*)^T \!\left[ H(z(s))\right]\!(z\!-\!z^*)ds 
\end{equation}
where $z=(x,y)$, $z^*=(x^*,y^*)$, $z(s)\!=\!z^* \!+\! s(z\!-\!z^*)$, and $H(z)$ is given by \eqref{eq:H}.
\item We then show that under the assumption \eqref{eq:condition-th16} $H(z)\preceq0$, and that for any 
$\tilde z = (\tilde\phi,\tilde P,\tilde \lambda,\tilde \nu_\G,\tilde \pi,\tilde\rho)\in\mathbb{R}^{2|\N|+3|\E|+|\G|+|\K|}$,
we have
\begin{align}\label{eq:null}
\tilde z^T H(z) \tilde z\!=\!0,\, \forall z\!\iff\! \tilde z\in \{\tilde z\in\mathbb{R}^\mathrm{Z}: \tilde\lambda\!=\!0,\tilde\nu_\G \!=\! 0, C_{\L}\tilde P\!=\!0 \}
\end{align}
where $\mathrm{Z}=2|\N|+3|\E|+|\G|+|\K|$.
\item We finally use \eqref{eq:null} and the invariance principle for Caratheodory systems~\cite{Bacciotti:2006uj} to show that
$\nu(t)\equiv 0$ and $\lambda(t)\equiv \lambda^*$.
\end{enumerate}
\vspace{.5ex}
The rest of the proof follows from steps 3 and 4 of Theorem \ref{thm:global-convergence}. 
\color{black}


We use $z=(x,y)$ and compactly express \eqref{eq:red-primal-dual-2} using
\begin{equation}\label{eq:dotz}
\dot z = Z[f(z)]^+_\rho
\end{equation}
where $Z=\blockdiag(X,Y)$ and 
\[
f(z) :=\left[
\begin{array}{c}
-\frac{\partial}{\partial x} L(x,y )^T\\
\frac{\partial}{\partial y} L(x,y )^T + g(x,y)
\end{array}
\right].
\]
Similarly, \eqref{eq:U} becomes $U(z)=\frac{1}{2} (z-z^*)^TZ^{-1}(z-z^*)$.

\vspace{1.5ex}
{\bf Proof of step 1:} We now recompute $\dot U(z)$ differenlty, i.e.
\begin{align}
&\dot U(z) \!=\!
\frac{1}{2}((z\!-\!z^*)^T[f(z)]^+_\rho \!+\! {[f(z)]^+_\rho}^T(z\!-\!z^*))\label{eq:dotUz-0}\\
&\!\leq\!\frac{1}{2}((z\!-\!z^*)^Tf(z) \!+\! f(z)^T(z\!-\!z^*))=(z\!-\!z^*)^Tf(z)\label{eq:dotUz-1}\\
&\!=\!\int_0^1\! (z\!-\!z^*)^T \!\left[\frac{\partial}{\partial z} f(z(s))\right]\!(z\!-\!z^*)ds \!+\!(z\!-\!z^*)^Tf(z^*)\label{eq:dotUz-1'}\\
&\!\leq\! \frac{1}{2}\int_0^1 \!(z\!-\!z^*)^T\! \left[\frac{\partial}{\partial z} f(z(s))^T \!+\! \frac{\partial}{\partial z} f(z(s))\right]\!(z\!-\!z^*)ds \label{eq:dotUz-2} \\
&\!=\! \int_0^1 \!(z\!-\!z^*)^T \!\left[ H(z(s))\right]\!(z\!-\!z^*)ds \label{eq:dotUz-integral-2}
\end{align}
where \eqref{eq:dotUz-0} follows from \eqref{eq:dotz}, \eqref{eq:dotUz-1} from \eqref{eq:projection-property}, and \eqref{eq:dotUz-1'} form 
the fact that $f(z)-f(z^*)=\int_0^1\ddx{z}f(z(s))(z-z^*)ds,$
where 
\begin{align*}
\frac{\partial}{\partial z} f(z) &=\left[
\begin{array}{cc}
-\frac{\partial^2}{\partial x^2} L(x,y ) & -\frac{\partial^2}{\partial x\partial y} L(x,y )\\
\frac{\partial^2}{\partial x\partial y} L(x,y )^T& \frac{\partial^2}{\partial y^2} L(x,y ) 
\end{array}
\right]\\&+
\left[
\begin{array}{cc}
0 & 0\\
 \frac{\partial}{\partial x}g(x,y) & \frac{\partial}{\partial y}g(x,y)
\end{array}
\right].
\end{align*}
Finally, \eqref{eq:dotUz-2} follows from the fact that either $f_i(z^*)=0$, or $(z_i-z_i^*)=z_i\geq0$ and $f_i(z^*)<0$, which implies $(z-z^*)^Tf(z^*)\leq 0$.

Therefore, $H(z)$ in \eqref{eq:dotUz-integral} can be expressed as
\begin{align*}
&H(z)=\frac{1}{2}\left[\frac{\partial}{\partial z} f(z)^T + \frac{\partial}{\partial z} f(z)\right]\\
&\;=\left[
\begin{array}{cc}
-\frac{\partial^2}{\partial x^2} L(x,y ) & 0\\
0& \frac{\partial^2}{\partial y^2} L(x,y ) 
\end{array}\right]\\
&\;+\left[
\begin{array}{cc}
0 & \frac{1}{2} \frac{\partial}{\partial x}g(x,y)^T\\
\frac{1}{2}\frac{\partial}{\partial x}g(x,y)  &\frac{1}{2}\left(\frac{\partial}{\partial y}g(x,y)^T+\frac{\partial}{\partial y}g(x,y)\right) 
\end{array}\right]
\end{align*}
which using lemmas \ref{lem:second-order-derivatives} and \ref{lem:partial-derivatives} becomes
\begin{equation}\label{eq:H}
H(z)=
\begin{blockarray}{c@{\hspace{5pt}}c@{\hspace{5pt}}c@{\hspace{5pt}}c@{\hspace{5pt}}c@{\hspace{5pt}}cl}
  &\matindex{\phi}  &\matindex{(P,\lambda_\L)}  & \matindex{(\lambda_\G,\nu_\G)} & \matindex{(\pi,\rho)} & &\\
\begin{block}{[c@{\hspace{5pt}}c@{\hspace{5pt}}c@{\hspace{5pt}}c@{\hspace{5pt}}c@{\hspace{5pt}}c]l}
& 0 & 0 & 0 & 0 & & \matindex{\phi}\\
& 0 & H_{P,\lambda_\L}(z) & 0 & 0 & & \matindex{(P,\lambda_\L)}\\
& 0 & 0 & H_{\lambda_\G,\nu_\G}(z) & 0 & & \matindex{(\lambda_\G,\nu_\G)}\\
& 0 & 0 & 0 & 0 & & \matindex{(\pi,\rho)} \\
\end{block}
\end{blockarray}
\end{equation}
where
\begin{align*}
&H_{P,\lambda_\L}(z)=\\
&
\begin{blockarray}{c@{\hspace{2.5pt}}c@{\hspace{2.5pt}}c@{\hspace{2.5pt}}c}
\begin{block}{[c@{\hspace{2.5pt}}c@{\hspace{2.5pt}}c@{\hspace{2.5pt}}c]}
&  -C_\L^T (D_\L+d'_\L)^{-1}C_\L 								& 	-\frac{1}{2}C_\L^T (D_\L+d'_\L)^{-1}\delta A_\L 	& \\
& 	-\frac{1}{2}\delta A_\L (D_\L+d'_\L)^{-1} C_\L 		&\; 	-(D_\L+\delta A_\L) (D_\L+d'_\L)^{-1}d'_\L 			\;& \\
\end{block}
\end{blockarray}\\
&\text{ and }\quad H_{\lambda_\G,\nu_\G}(z)=
\begin{blockarray}{c@{\hspace{5pt}}c@{\hspace{5pt}}c@{\hspace{5pt}}c}
\begin{block}{[c@{\hspace{5pt}}c@{\hspace{5pt}}c@{\hspace{5pt}}c]}
& -d'_\G 							& \frac{1}{2}\delta A_\G-d'_\G	& \\
&  \frac{1}{2}\delta	 A_\G-d'_\G 		& -(d'_\G+D_\G) 									& \\
\end{block}
\end{blockarray}.
\end{align*}

It will prove useful in the next step to rewrite $H_{P,\lambda_L}(z)$ using
\begin{equation}\label{eq:th16-1}
H_{P,\lambda_L}(z)=\tilde C^T \tilde D^{\frac{1}{2}}(z) 
\tilde H(z) \tilde D^{\frac{1}{2}}(z)\tilde C
\end{equation}
where  $$\tilde C\!=\!\blockdiag(C_\L,I),\,\tilde D(z) = \blockdiag(D_\L+d'_\L, D_\L+d'_\L)^{-1}$$
\begin{align*}
\text{ and }\quad
\tilde H(z)&=\left[
\begin{array}{cc}
 -I & -\frac{1}{2}\delta A_\L\\
 -\frac{1}{2}\delta A_\L  & -(D_\L+\delta A_\L)d'_\L
\end{array}
\right].\qquad\qquad\qquad\qquad
\end{align*}
Notice that since $\tilde D(z) \succ 0$, $\tilde D^{\frac{1}{2}}(z)$  in \eqref{eq:th16-1} always exists.

\vspace{1.5ex}
{\bf Proof of step 2:}
To show that $H(z)\preceq 0 $ and \eqref{eq:null} holds, it is enough to show that 
\begin{equation}\label{eq:condition-step-2}
\tilde H(z)\prec0 \quad\text{ and }\quad H_{\lambda_\G,\nu_\G}(z)\prec 0,  \quad\forall z.
\end{equation}

To see this, assume for now that \eqref{eq:condition-step-2} holds.
Then, using \eqref{eq:th16-1} it follows that $H_{P,\lambda_\L}(z)\preceq0$, which implies by \eqref{eq:H} and $H_{\lambda_\G,\nu_\G}(z)\prec 0$ that $H(z)\preceq0$.
Moreover, $\tilde z^T H(z) \tilde z=0$ $\forall z$ if and only if 
\begin{align}\label{eq:th16-2}
[\,\tilde P^T \, \tilde\lambda_\L^T\,] H_{P,\lambda_\L}(z) [\,\tilde P^T \, \tilde\lambda_\L^T\,]^T=0
\end{align}
and
\begin{align}\label{eq:th16-3}
[\,\tilde\lambda_\G^T\,\tilde \nu_\G^T \,] H_{\lambda_\G,\nu_\G}(z) [\,\tilde\lambda_\G^T\,\tilde \nu_\G^T \,]^T=0.
\end{align}
Therefore using \eqref{eq:th16-1} it follows that \eqref{eq:th16-2} and $\tilde H(z)\prec 0$ $\forall z$ implies that $C_\L\tilde P=0$ and $\tilde\lambda_\L=0$. Similarly, $H_{\lambda_\G,\nu_\G}(z)\prec 0$ $\forall z$ and \eqref{eq:th16-3} implies $\tilde\lambda_\G=\tilde\nu_\G=0$. This finishes the proof of \eqref{eq:null}.
It remains to show that \eqref{eq:condition-step-2} holds whenever \eqref{eq:condition-th16} holds.\\[1.5ex]
\noindent
{\it Proof of $\tilde H(z)\prec 0$:}
By definition of negative definite matrices, $\tilde H(z)\prec 0$ if and only if all the roots of the characteristic polynomials
\begin{align*}
&p_i(\mu_i)=(\mu_i+1)(\mu_i+(D_i+\delta a_i)d'_i)-{\delta a_i^2}/{4}\\
&= \mu_i^2+ (1+ (D_i+\delta a_i)d'_i)\mu_i+ (D_i+\delta a_i)d'_i-{\delta a_i^2}/{4}
\end{align*}
are negative for every $i\in \L$ and $\forall z$ (recall $d'_i$ depends on $z$). 

Thus, applying Ruth-Hurwitz stability criterion we get the following necessary and sufficient condition:
\begin{subequations}\label{eq:th14-1}
\begin{align}
\delta a_i^2-4(D_i+\delta a_i) d'_i<0\label{eq:th14-1-a}\\
 1 + (D_i+\delta a_i)d'_i > 0\label{eq:th14-1-b}
\end{align}
\end{subequations}
for every $i\in\L$. 

Now, equation \eqref{eq:th14-1-a} can be equivalently rewritten as:
\begin{equation}\label{eq:th14-2}
2(d'_i -\sqrt{d'_i(d'_i+D_i)}) < \delta a_i < 2(d'_i +\sqrt{d'_i(d'_i+D_i)}).
\end{equation}
Since $d'_i\in[\underline{d}',\bar{d}']$, $D_i\geq D_{\min}$ and the function 
$ 
x-\sqrt{x(x+y)}
$
 is decreasing in both $x$ and $y$ for $x,y\geq 0$, then
\begin{align*}
2(d'_i -\sqrt{d'_i(d'_i+D_i)})&\leq 2(\underline d' -\sqrt{\underline d'(\underline d'+D_{\min})}).
\end{align*}
Similarly, since $x+\sqrt{x(x+y)}$ is increasing for $x,y\geq0$, 
\begin{align*}
2(d'_i +\sqrt{d'_i(d'_i+D_i)})&\geq 2(\underline{d}'+\sqrt{\underline{d}'(\underline{d}'+D_{\min})}).
\end{align*}
Therefore, \eqref{eq:th14-1-a} holds whenever $\delta a_i$ satisifes \eqref{eq:condition-th16}.

Finally,  \eqref{eq:th14-1-b} holds whenever $\delta a_i  > - \frac{1}{d_i'} - D_i$ which in particular holds if  
$
\delta a_i  > -  D_{\min}.
$
The following calculation shows that $2(\underline d' -\sqrt{\underline d'(\underline d'+D_{\min})})>-D_{\min}$ 
which implies that \eqref{eq:th14-1-b} holds under condition \eqref{eq:condition-th16}:
\begin{align*}
2(\underline d' -\sqrt{\underline d'(\underline d'+D_{\min})}) 	& >-D_{\min} 							&\iff&\\
  \sqrt{\underline d'(\underline d'+D_{\min})}				& <\underline d' +\frac{D_{\min}}{2} 			&\iff&\\[-.5ex]
  \underline d'(\underline d'+D_{\min}) 				< {\underline d'}^2 +&\frac{D_{\min}^2}{4} +\underline d'D_{\min} &\iff&
\;\;\;  0 < \frac{D_{\min}^2}{4}.
\end{align*}
Therefore \eqref{eq:th14-1} holds  whenever \eqref{eq:condition-th16} holds.\\[1.5ex]
\noindent
{\it Proof of $H_{\nu_\G,\lambda_\G}(z)\prec0$:}
Similarly, we can show that {\it all} the eigenvalues of $H_{\nu_\G,\lambda_\G}(z)$ are the roots of the polynomials
\begin{align*}
p_i(\mu_i)&=(\mu_i+ D_i + d'_i)(\mu_i +d'_i) - (\frac{\delta a_i}{2}-d'_i)^2\\[-.5ex]
&=\mu_i^2 + (D_i + 2d'_i)\mu_i + (D_i+\delta a_i)d'_i - \frac{\delta a_i^2}{4}
\end{align*}
which, since $D_i+2d'_i>0$, are negative if and only if \eqref{eq:th14-1-a} is satisfied $\forall i\in \G$. Therefore,
\eqref{eq:condition-th16} also guarantees that $H_{\nu_\G,\lambda_\G}\prec 0$.

\vspace{1.5ex}
{\bf Proof of step 3:}  Since by Step 2 $H(z)\preceq 0$ $\forall z$, \eqref{eq:dotUz-integral} implies that $\dot U\leq 0$ whenever \eqref{eq:condition-th16} holds.
 Thus, we are left to apply again the  invariance principle for Caratheodory systems~\cite{Bacciotti:2006uj} and characterize its invariant set $M$~\eqref{eq:M}.

Let $\delta z = (z(t)-z^*)$ and similarly define $\delta P=(P(t) - P^*)$,  $\delta \lambda_\L=\lambda_\L(t) - \lambda_\L^*$,
$\delta \lambda_\G=\lambda_\G(t) - \lambda_\G^*$ and $\delta \nu_\G=\nu_\G(t) - \nu_\G^*$.
Then since $\dot U\equiv 0$ iff $\delta z^T H(z)\delta z\equiv0$,
then it follows from \eqref{eq:null} that $z(t)\in M$ if and only if $C_\L \delta P\equiv 0 $, $\delta \lambda\equiv 0$ and $\delta \nu_\G \equiv 0$.

This implies that $C_\L P(t)\equiv C_\L P^*$, $\lambda(t)\equiv \lambda^*$ and $\nu_\G(t)\equiv\nu_\G^*=0$, which in particular also implies that $\nu^*_\L(x(t),y(t))=\nu^*_\L(C_\L P(t),\lambda(t))\equiv\nu^*_\L(C_\L P^*,\lambda^*)=0$.
Therefore we have shown that $z(t)\in M$ if and only if $\lambda(t)\equiv\lambda^*$ and $\nu(t)\equiv 0$ which finalizes  Step 3.\\[1.5ex]
As mentioned before, the rest of the proof follows from {\bf steps 2 and 3} of Theorem \ref{thm:global-convergence}.
\end{IEEEproof}
\color{black}

\color{black}
\section{Framework Extensions}\label{sec:framework-extensions}

In this section we extend the proposed framework to derive controllers that enhance the solution described before. More precisely, we will show how we can modify our controllers in order to account for buses that have zero power injection (Section \ref{ssec:zero-injection-buses}) and how to  fully distribute the implementation of the inter-area flow contraints (Section \ref{ssec:inter-area-flow-constraints}).

\subsection{Zero Power Injection Buses}\label{ssec:zero-injection-buses}
We now show how our design framework can be extended to include buses with zero power injection.
Let $\Z$ be the set of buses that have neither generators nor loads. Thus,  we consider a power network whose dynamics are described by 
\begin{subequations}\label{eq:swing-theta-Z}
\begin{align}
\!\!\dot \theta_{\G\cup\L} &\!=\! \omega_{\G\cup\L}  										\label{eq:sw-theta-Z-a} \\
\!\!M_\G \dot\omega_\G & \!=\! P^{in}_\G \!\!-\!(d_\G \!+\! D_\G\omega_\G) \!-\!L_{B,(\G,\N)}\theta\ \label{eq:sw-theta-Z-b}\\
\!\!0 & \!=\! P^{in}_\L\!\!-\!(d_\L \!+\! D_\L\omega_\L) \!-\!L_{B,(\L,\N)}\theta		\label{eq:sw-theta-Z-c}\\
\!\!0 & \!=\! -\! L_{B,(\Z,\N)}\theta									\label{eq:sw-theta-Z-d}
\end{align}
\end{subequations}
where $L_{B,(S,S')}$ is the sub-matrix of $L_B$ consisting of the rows in $S$ and columns in $S'$.

We will use Kron reduction to eliminate \eqref{eq:sw-theta-Z-d}. Equation \eqref{eq:sw-theta-Z-d} implies that the ($\theta_i$, $i\in\Z$) is uniquely determined by the buses adjacent to $\Z$, i.e. $\theta_\Z =  L_{B,(\Z,\Z)}^{-1}L_{B,(\Z,\G\cup\L)}\theta_{\G\cup\L}$.
Thus we can rewrite \eqref{eq:swing-theta-Z} using only $\theta_{\G\cup\L}$ which gives
\begin{subequations}\label{eq:swing-theta-sharp}
\begin{align}
\!\!\dot \theta_{\G\cup\L} &\!=\! \omega_{\G\cup\L}  										\label{eq:sw-theta-sharp-a} \\
\!\!M_\G \dot\omega_\G & \!=\! P^{in}_\G \!\!-\!(d_\G \!+\! D_\G\omega_\G) \!-\!L^\sharp_{B,(\G,\G\cup\L)}\theta_{\G\cup\L} \label{eq:sw-theta-sharp-b}\\
\!\!0 & \!=\! P^{in}_\L\!\!-\!(d_\L \!+\! D_\L\omega_\L) \!-\!L^\sharp_{B,(\L,\G\cup\L)}\theta_{\G\cup\L}		\label{eq:sw-theta-sharp-c}
\end{align}
\end{subequations}
where 
$
L^\sharp_B = L_{B,(\G\cup\L,\G\cup\L)} - L_{B,(\G\cup\L,\Z)} L_{B,(\Z,\Z)}^{-1}L_{B,(\Z,\G\cup\L)}
$
 is the Schur complement of $L_B$ after removing the rows and columns corresponding to $\Z$. The matrix $L^\sharp_B$ is also a Laplacian of a reduced graph $G({\G\cup\L},\E^\sharp)$ and therefore it can be expressed as  
$
L^\sharp_B = C^{\sharp T}B^\sharp C^{\sharp T}
$
 where $C^\sharp$ is the incidence matrix of $G({\G\cup\L},\E^\sharp)$  and $B^\sharp=\diag(B^\sharp_{ij})_{ij\in\E^\sharp}$ are the line susceptances of the reduced network.

This reduction allows to use \eqref{eq:swing-theta-sharp} (which is equivalent to \eqref{eq:swing-theta}) to also model networks that contain buses with zero power injection. The only caveat is that some of line flows of the vector $BC^T\theta$ are no longer present in $B^\sharp C^{\sharp T}\theta_{\G\cup\L}$ -- when a bus is eliminated using Kron reduction, its adjacent lines $B_{e}$, $e\in\E$, are substituted by an equivalent clique with new line impedances $B^\sharp_{e'}$, $e'\in\E^\sharp$.
As a result, some of the constraints \eqref{eq:welfare-d}-\eqref{eq:welfare-e} would no longer have a physical meaning if we directly substitute $BC^T\theta$ with $B^\sharp C^{\sharp T}\theta_{\G\cup\L}$ in \eqref{eq:swing-theta}. 

We overcome this issue by showing that each original $B_{ij}(\theta_i-\theta_j)$ in $G(\N,\E)$ can be replaced by a linear combination of line flows $B^\sharp_{i'j'}(\theta_{i'}-\theta_{j'})$ of the reduced network $G(\G\cup\L,\E^\sharp)$.

For any $\theta$ satisfying \eqref{eq:sw-theta-Z-d} we have
\begin{align*}
L_B\theta=\left[\begin{array}{c} q_{\G\cup\L}\\ \0_{|\Z|} \end{array}\right]=\left[\begin{array}{c}L_B^\sharp\theta_{\G\cup\L}\\ \0_{|\Z|} \end{array}\right].
\end{align*}
Thus it follows that 
\begin{align}
&BC^T\theta = BC^T L_B^\dagger \left[\begin{array}{c} q_{\G\cup\L}\\ \0_{|\Z|} \end{array}\right]=BC^T L_{B,(\N,\G\cup\L)}^\dagger q_{\G\cup\L}\nonumber\\
&\!=\!BC^T L_{B,(\N,\G\cup\L)}^\dagger  C^\sharp B^\sharp C^{\sharp T} \theta_{\G\cup\L} \!:=\! A^\sharp B^\sharp C^{\sharp T} \theta_{\G\cup\L}\label{eq:Asharp}
\end{align}
 where $L_B^\dagger$ is the pseudo-inverse of $L_B$.
 
Therefore, by substituting $BC^T\theta$ with $A^\sharp B^\sharp C^{\sharp T}\theta_{\G\cup\L}$ in \eqref{eq:welfare} and repeating the procedure of Section \ref{sec:frequency-preserving-olc} we obtained a modified version of \eqref{eq:load-control} in which \eqref{eq:load-control-a}-\eqref{eq:load-control-e} becomes
\begin{subequations}\label{eq:load-control-3}
\begin{align}
\dot\lambda_i &=\zeta_i^\lambda\Big(P^m_i -d_i -\sum_{j\in\N^\sharp_i}B^\sharp_{ij}(\phi_i-\phi_j)\Big)		\label{eq:load-control-3-a}\\
\dot\pi_k &= \zeta^\pi_k\Big(\!\!\!\!\!\!\!\!\!\!\!\!\sum_{\;\;\;\;\;\;e\in\B_k,ij\in\E^\sharp}\!\!\!\!\!\!\!\!\!\!\!\!\hat C_{k,e}A^\sharp_{e,ij}B^\sharp_{ij}(\phi_i-\phi_j)-\hat P_{k}\Big)		\label{eq:load-control-3-b}\\
\dot\rho_{e}^+ &= \zeta_{e}^{\rho^+}\Big[\sum_{ij\in\E^\sharp} A_{e,ij}^\sharp B^\sharp_{ij}(\phi_i-\phi_j) - \overline{P_{e}}\Big]^+_{\rho_{e}^+}							\label{eq:load-control-3-c}\\
\dot\rho_{e}^- &= \zeta_e^{\rho^-}\Big[\underline {P_{e}} -\sum_{ij\in\E^\sharp}A_{e,ij}^\sharp B^\sharp_{ij}(\phi_i-\phi_j) \Big]^+_{\rho_e^-}							\label{eq:load-control-3-d}\\
\dot \phi_i&= \chi_i^\phi\Big(\!\sum_{j\in\N^\sharp_i} B^\sharp_{ij}(\lambda_i\!-\!\lambda_j)\!-\!\sum_{e\in\E}C^\sharp_{i,e}B^\sharp_e \!\!\!\!\!\sum_{e'\in\E,\,k\in\K}\!\!\!\!\!A^\sharp_{e,e'}\hat C_{k,e'}\pi_k 		\nonumber\\
				& \qquad\quad-\!\sum_{e\in\E}C^\sharp_{i,e}B^\sharp_e \sum_{e'\in\E}A^\sharp_{e,e'}(\rho_{e'}^{+}-\rho_{e'}^{-}) \Big) 																					\label{eq:load-control-3-e}
\end{align}
\end{subequations}
where \eqref{eq:load-control-3-a} and \eqref{eq:load-control-3-e}  are for $i\in\G\cup\L$, \eqref{eq:load-control-3-b} is for $k\in\K$, and \eqref{eq:load-control-3-c} and \eqref{eq:load-control-3-d} are for the original lines $e\in \E$.

\begin{figure}[htp]\centering
\includegraphics[width=.425\columnwidth]{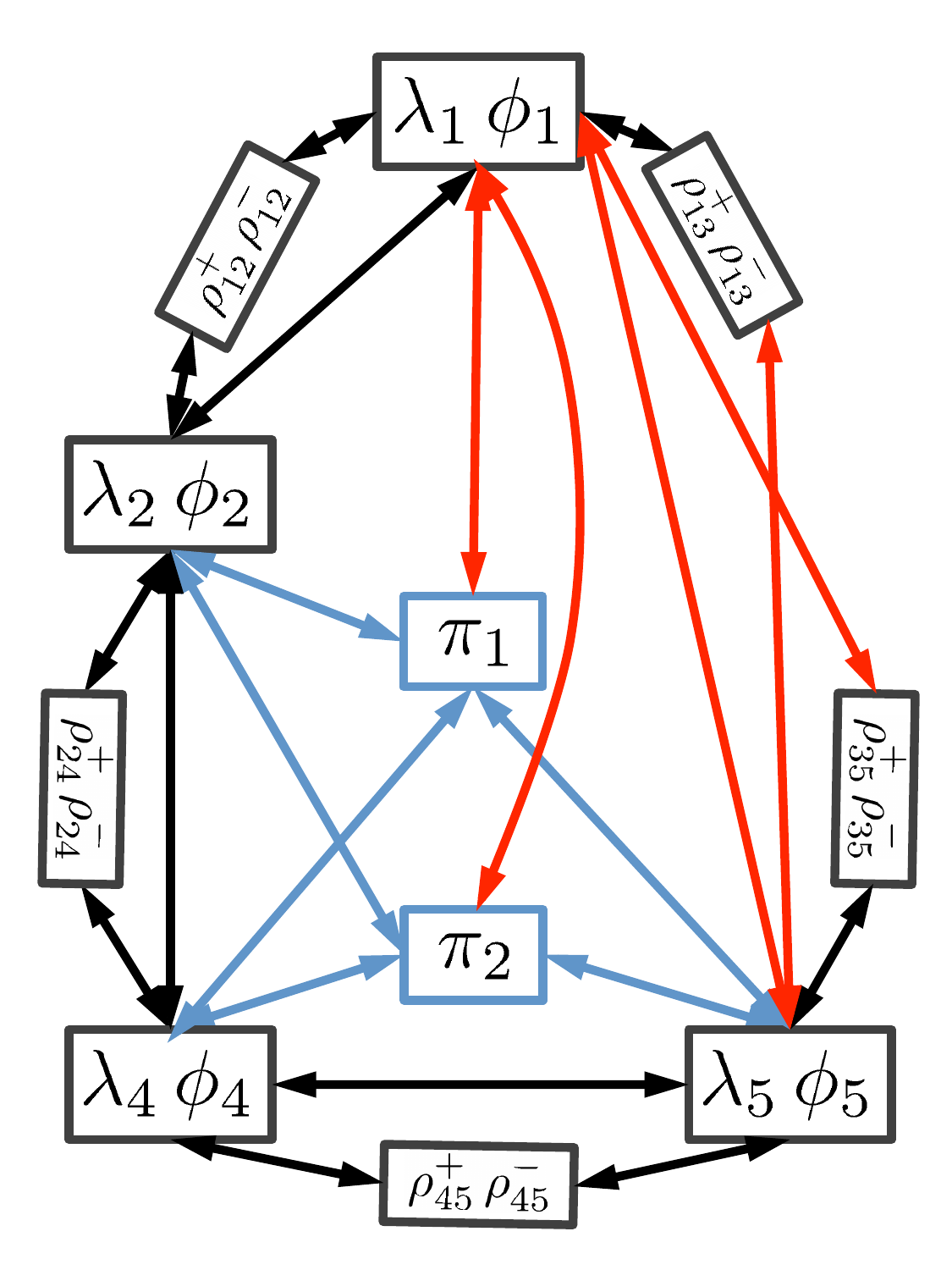}
\includegraphics[width=.425\columnwidth]{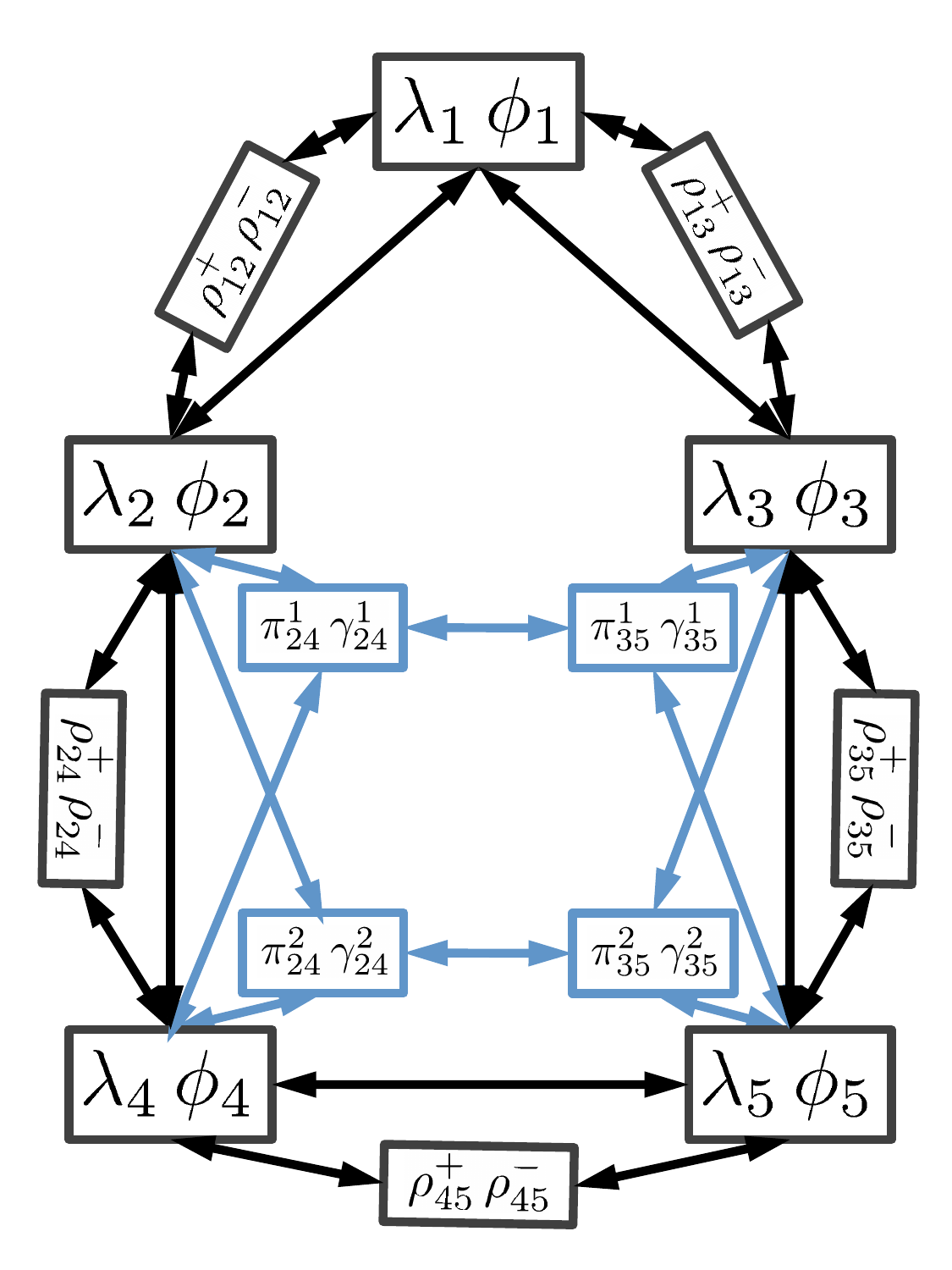}
\caption{Communication requirements for the power network in Fig. \ref{fig:diagram} when bus $3$ has no injection as described in Section \ref{ssec:zero-injection-buses} (left), and for the distributed inter-area flow constraint formulation of Section \ref{ssec:inter-area-flow-constraints} (right). }\label{fig:diagram-reformulation}
\end{figure}

It can be shown that the analysis described in Sections \ref{sec:analysis} and \ref{sec:implementation-details} still holds under this extension.
 \begin{remark}
The only additional overhead incurred by the proposed extension is the need for communication between buses that are adjacent on the graph $G(\G\cup\L,\E^\sharp)$ and were not adjacent in $G(\N,\E)$ (see Figure \ref{fig:diagram-reformulation} for an illustration).
\end{remark}

\subsection{Distributed Inter-area Flow Constraints}\label{ssec:inter-area-flow-constraints}

We now show how we can fully distribute the implementation of the inter-area flow constraints. The procedure is analogous to Section \ref{ssec:zero-injection-buses} and therefore we will only limit to describe what are the modifications that need to be done to~\eqref{eq:welfare} in order to obtain controllers that are fully distributed.

We define for each area $k$ an additional graph $G(\B_k,\E^k)$ where we associate each boundary edge $e\in\B_k$ with a node and define \emph{undirected} edges $\{e,e'\}\in \E^k$ that describe the communication links between $e$ and $e'$. Using this formulation, we decompose equation \eqref{eq:inter-area-constraint-2} for each $k$ into $|\B_k|$ equations
\begin{align}\label{eq:decomposed-iac}
\hat C_{k,e}P_e -\frac{\hat P_k}{|\B_k|}=\!\!\! \sum_{e':\{e,e'\}\in\B_k} \!\!\!(\gamma_e-\gamma_{e'}),\quad&e\in\B_k,\,k\in \K
\end{align} 
where $\gamma_e$ is a new primal variable that aims to guarantee indirectly \eqref{eq:inter-area-constraint-2}. In fact, it is easy to see by summing \eqref{eq:decomposed-iac} over $e\in\B_k$ that
\begin{align*}
\sum_{e\in\B_k}\left(\hat C_{k,e}P_e -\frac{\hat P_k}{|\B_k|}\right)=\sum_{e\in\B_k}\sum_{e':\{e,e'\}\in\B_k} \!\!\!(\gamma^k_e-\gamma^k_{e'})=0.
\end{align*}
which is equal to \eqref{eq:inter-area-constraint-2}.

Therefore, since whenever \eqref{eq:inter-area-constraint-2} holds, one can find a set of $\gamma_e$ satisfying \eqref{eq:decomposed-iac}, then we can substitute \eqref{eq:welfare-d} with \eqref{eq:decomposed-iac}.
If we let $\pi^k_{ij}$ be the Lagrange multiplier associated with \eqref{eq:decomposed-iac}, then by replacing \eqref{eq:load-control-b} and \eqref{eq:load-control-e} with
\begin{subequations}\label{eq:load-control-4}
\begin{align}
\dot \pi^k_{ij} \!&=\! \zeta_{k,ij}^\pi\!\Big(\!\hat C_{k,ij}B_{ij}(\phi_i\!-\!\phi_j) \!-\!\frac{\hat P_k}{|\B_k|}\!-\!\!\!\!\!\!\!\!\!\!\!\!\!\!\sum_{\;\;\;\;\;e:\{ij,e\}\in\B_k} \!\!\!\!\!\!\!\!\!\!\!\!\!\gamma^k_{ij}\!-\!\gamma^k_{e}\!\Big)\label{eq:load-control-4-a}\\
\dot \gamma^k_{ij}\!&=\! \chi_{k,ij}^\gamma\Big ( \!\!\!\!\!\!\!\!\!\!\!\!\!\!\!\!\sum_{\quad\;\;\;\;e:\{ij,e\}\in\B_k} \!\!\!\!\!\!\!\!\!\!\!\!\!\pi^k_{ij}\!-\!\pi^k_{e}\;\Big)\label{eq:load-control-4-b}\\
\dot \phi_i&= \chi_i^\phi\Big(\sum_{j\in\N_i}\!\! B_{ij}(\lambda_i\!-\!\lambda_j)\!-\! \sum_{k\in\K,\,e\in\B_k}\!\!C_{i,e}B_e\hat C_{k,e}\pi_e^k 		\nonumber\\
				& \qquad\quad- \sum_{e\in\E}C_{i,e}B_e(\rho_e^{+}-\rho_e^{-}) \Big) 																					\label{eq:load-control-4-c}
\end{align}
\end{subequations}
we can distribute the implementation of the inter-area flow constraint. Figure \ref{fig:diagram-reformulation} shows how the communication requirements are modified by this change.
If each boundary bus has only one incident boundary edge, i.e. if  $\sum_{k\in\K,\,e\in\B_k}\!\!C_{i,e}B_e\hat C_{k,e}\pi_e^k$ has at most one term, the convergence results of sections \ref{sec:analysis} and \ref{sec:implementation-details} extend to this case.

\color{black}

\section{Numerical Illustrations}
\label{sec:numerical-illustrations}
\begin{figure}[tp]
\centering
\includegraphics[width=.95\columnwidth]{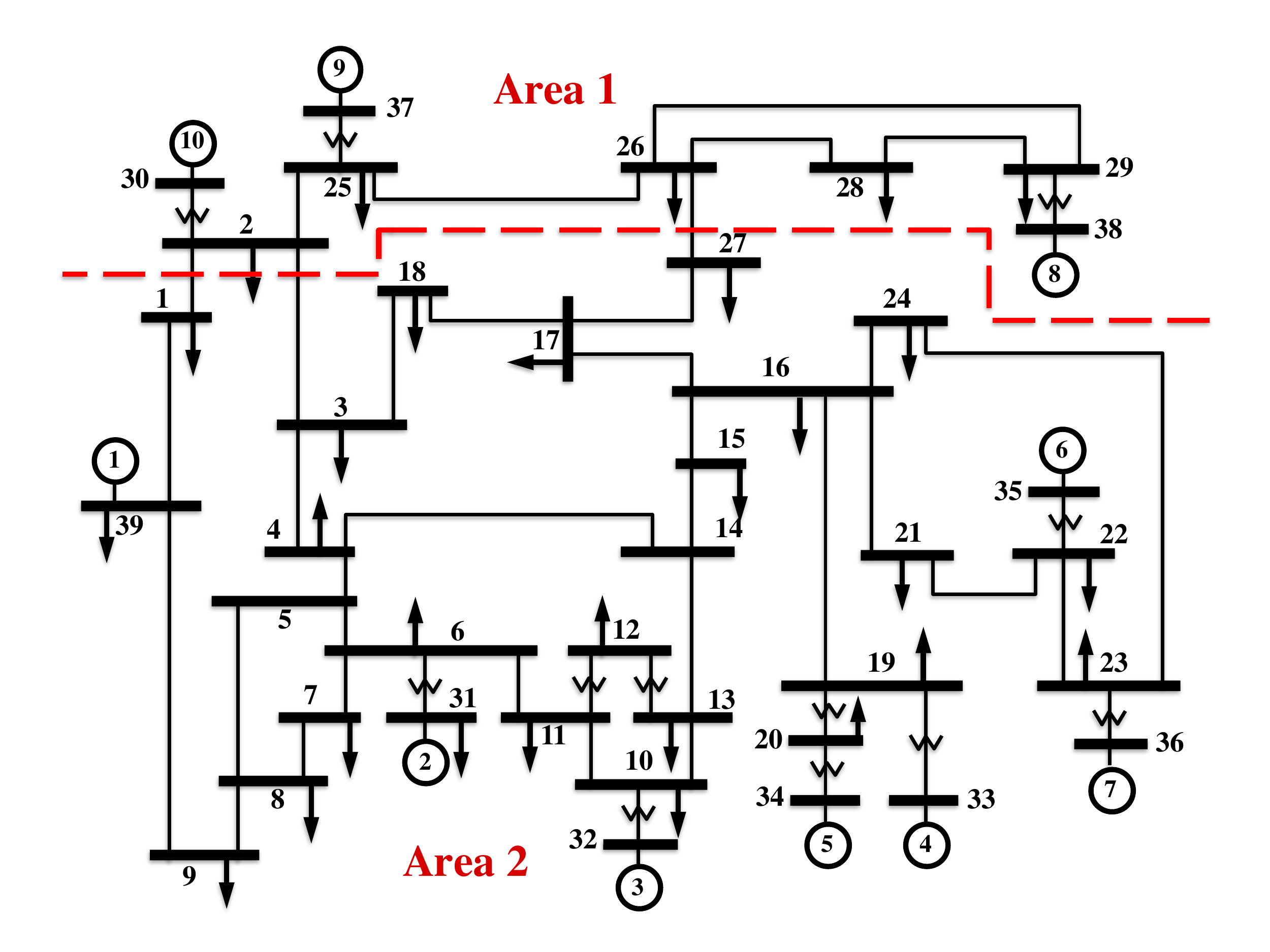}
\caption{IEEE 39 bus system: New England}\label{fig:new_england}
\end{figure}

We now illustrate the behavior of our control scheme. We consider the widely used IEEE 39 bus system, shown in Figure \ref{fig:new_england}, to test our scheme. We assume that the system has two independent control areas that are connected through lines $(1,2)$, $(2,3)$ and $(26,27)$.
The network parameters as well as the initial stationary point  (pre fault state) were obtained from the Power System Toolbox~\cite{chow1992toolbox} data set. 

\begin{figure}[ht!]
\centering
\includegraphics[width=.8\columnwidth]{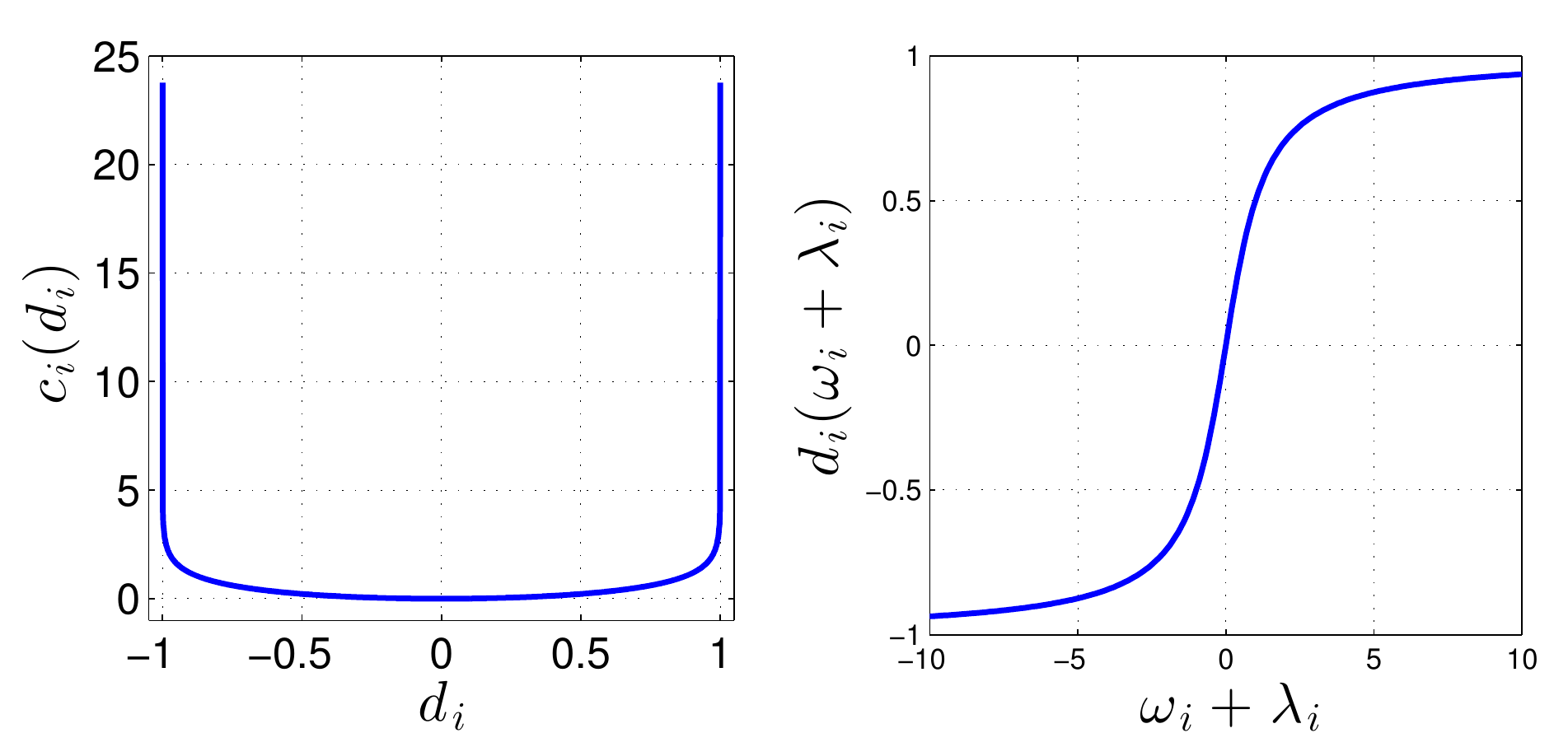}
\caption{Disutility $c_i(d_i)$ and load function $d_i(\omega_i+\lambda_i)$}\label{fig:ci_di}
\end{figure}

Each bus is assumed to have a controllable load with $\D_i=[-d_{\max},d_{\max}]$, with $d_{\max}=1$p.u. on a  $100$MVA base
with $c_i(\cdot)$ and $d_i(\cdot)$ as shown in Figure \ref{fig:ci_di}.

Throughout the simulations we assume that the aggregate generator damping and load frequency sensitivity parameter $D_i=0.2$ $\forall i\in\N$ and $\chi_i^\phi=\zeta^\lambda_i=\zeta^\pi_k=\zeta^{\rho^+}_e=\zeta^{\rho^-}_e=1$, for all $i\in\N$, $k\in\K$ and $e\in\E$. These parameter values do not affect convergence, but in general they will affect the convergence rate. 
The values of $P^{in}$ are corrected so that they initially add up to zero by evenly distributing the mismatch among the load buses.
$\hat P$ is obtained from the starting stationary condition. We initially set $\overline P$ and $\underline P$ so that they are not binding. 

\begin{figure}[h!]
\includegraphics[width=\columnwidth]{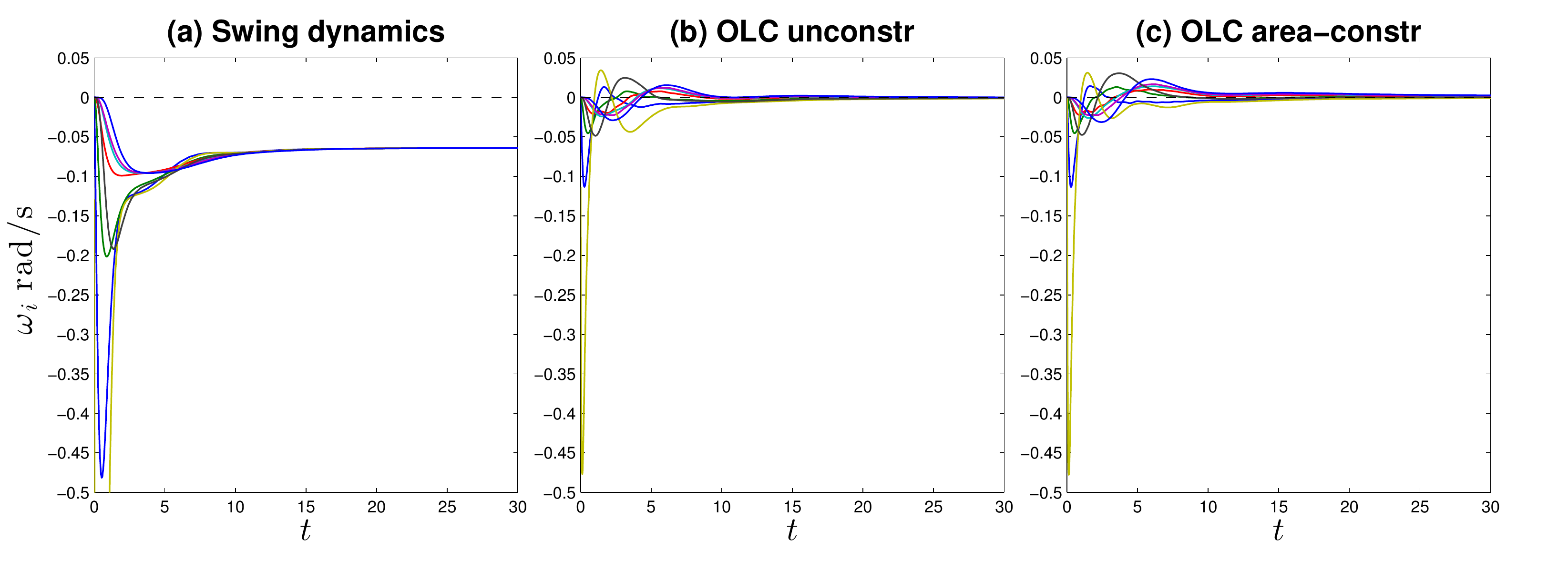}
\caption{Frequency evolution: Area 1}\label{fig:omegai-a1}
\end{figure}
\begin{figure}[h!]
\includegraphics[width=\columnwidth]{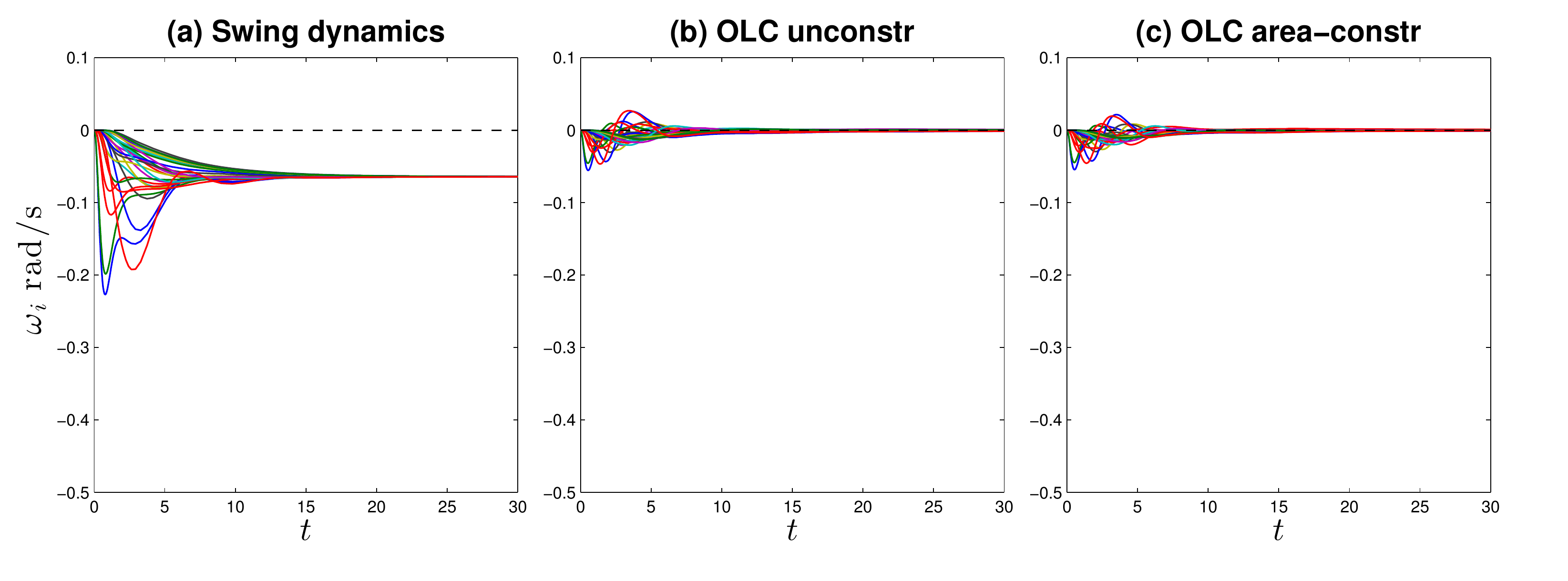}
\caption{Frequency evolution: Area 2}\label{fig:omegai-a2}
\end{figure}

We simulate the OLC-system as well as the swing dynamics~\eqref{eq:1} without load control ($d_i=0$), after introducing a perturbation at bus $29$ of $P^{in}_{29}=-2$p.u.. \added{In some scenarios we disable a few of the OLC constraints. This is achieved by fixing the corresponding Lagrange multiplier to be zero.}

Figures \ref{fig:omegai-a1} and \ref{fig:omegai-a2} show the evolution of the bus frequencies for the uncontrolled swing dynamics (a), the OLC system without inter-area constraints (b), and the OLC with area constraints (c). 
It can be seen that while the swing dynamics alone fail to recover the nominal frequency, the OLC controllers can jointly rebalance the power as well as recovering the nominal frequency. The convergence of OLC seems to be similar or even better than the swing dynamics, as shown in Figures \ref{fig:omegai-a1} and \ref{fig:omegai-a2}.

\begin{figure}[h!]\centering
\includegraphics[width=\columnwidth]{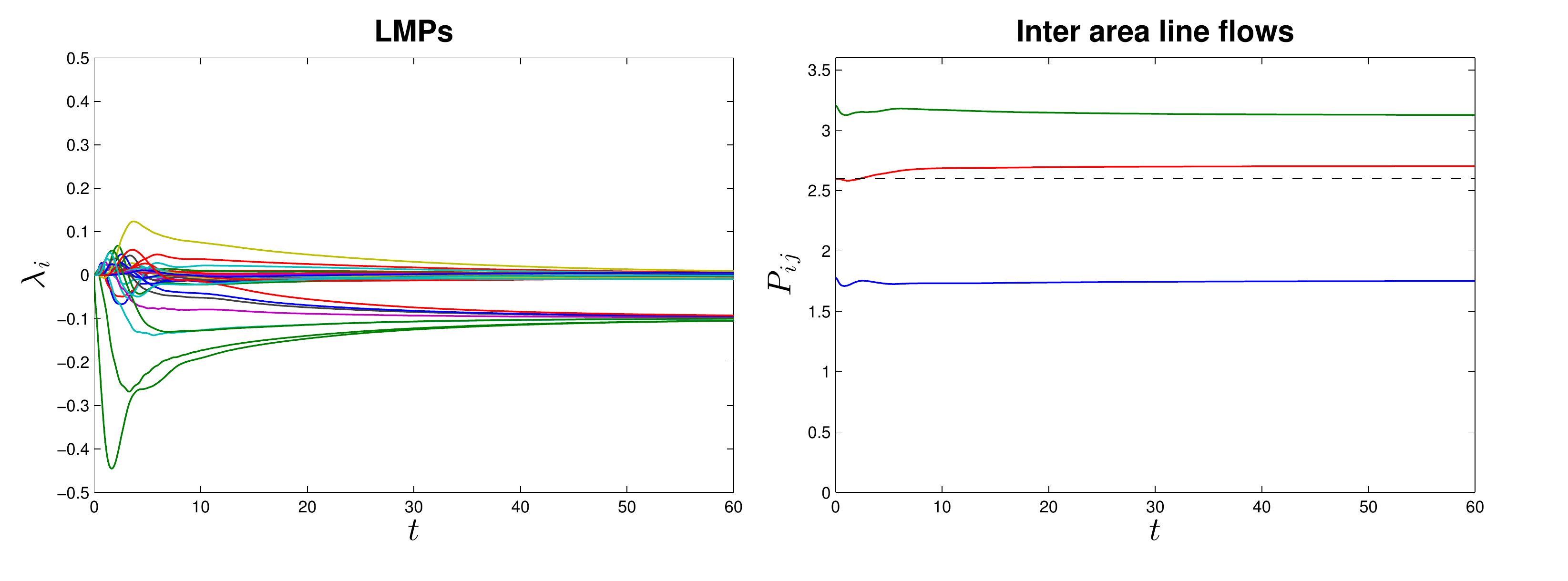}
\caption{LMPs and inter area line flows: no thermal limits}\label{fig:no-line-constr}
\end{figure}
\begin{figure}[h!]\centering
\includegraphics[width=\columnwidth]{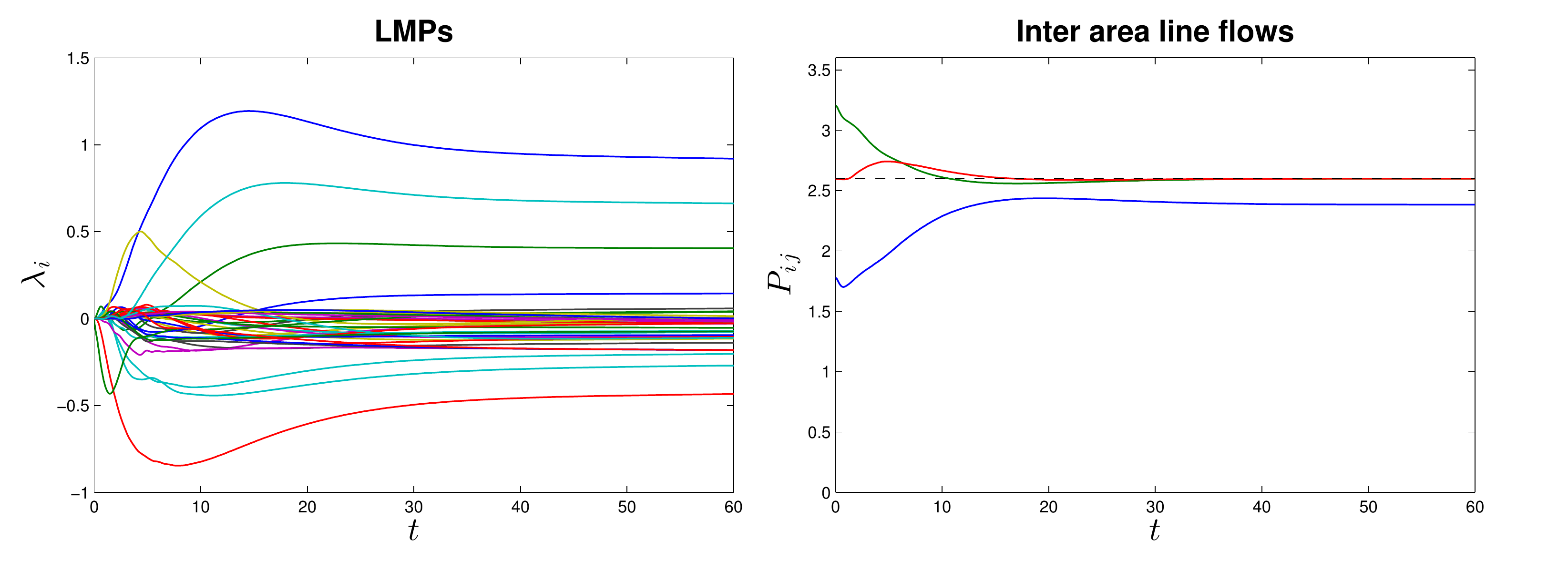}
\caption{LMPs and inter area line flows: with thermal limits}\label{fig:line-constr}
\end{figure}

Now, we illustrate the action of the thermal constraints by adding a constraint of $\overline P_e=2.6$p.u. and $\underline P_e=-2.6$p.u. to the tie lines between areas. Figure \ref{fig:no-line-constr} shows the values of the multipliers $\lambda_i$, that correspond to the Locational Marginal Prices (LMPs), and the line flows of the tie lines for the same scenario displayed in Figures \ref{fig:omegai-a1}c and \ref{fig:omegai-a2}c, i.e. without thermal limits. It can be seen that neither the initial condition, nor the new steady state satisfy the thermal limit (shown by a dashed line). However, once we add thermal limits to our OLC scheme, we can see in Figure \ref{fig:line-constr} that the system converges to a new operating point that satisfies our constraints.

\ifjournal
Finally, we show the conservativeness of the bound obtained in Theorem \ref{thm:global-convergence-2}. We simulate the  perturbed system \eqref{eq:1}, \eqref{eq:load-control-2-a} and \eqref{eq:load-control-b}-\eqref{eq:load-control-f} with under the same conditions as in Figure \ref{fig:no-line-constr}.
We set the scalars $a_i$s such that the corresponding $\delta a_i$s are homogeneous for every bus $i$. We also do not impose the bounds \eqref{eq:d'-bounds} on $d_i(\cdot)$ and use instead $d_i$ as described in Figure \ref{fig:ci_di}. This last assumption actually implies that the interval in \eqref{eq:condition-th16} is empty (because $\underline{d}'=0$).

\begin{figure}[htp]\centering
\includegraphics[width=\columnwidth]{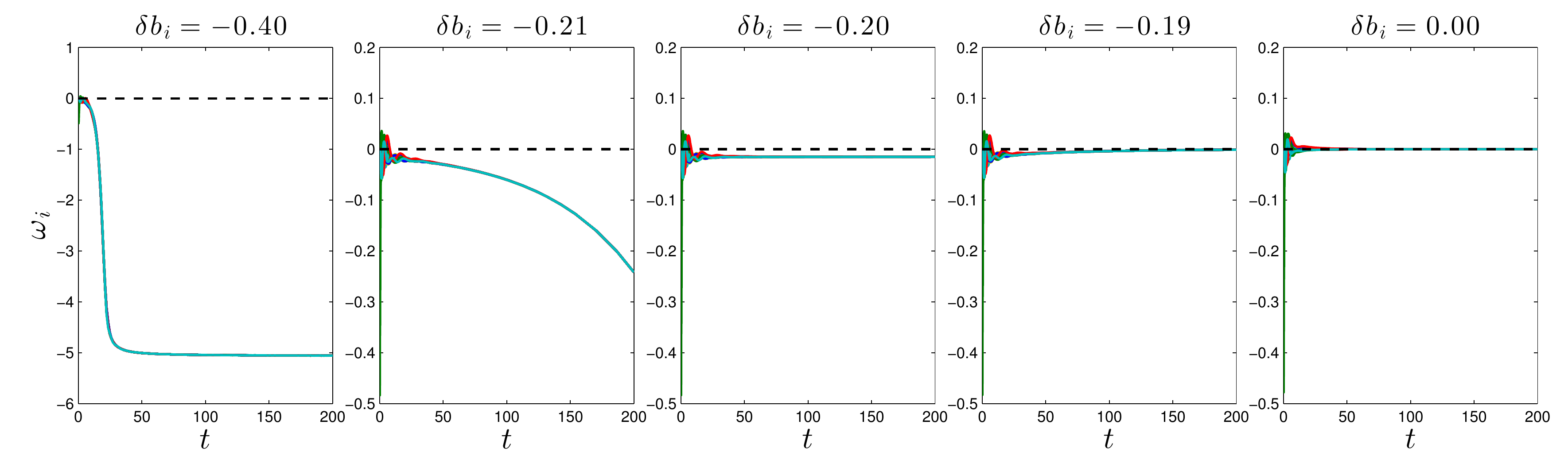}
\caption{Bus frequency evolution for homogeneous perturbation $\delta a_i\in\{-0.4,-0.21,-0.2,-0.19,0.0\}$}\label{fig:perturbed-system-1}
\end{figure}
\begin{figure}[htp]\centering
\includegraphics[width=\columnwidth]{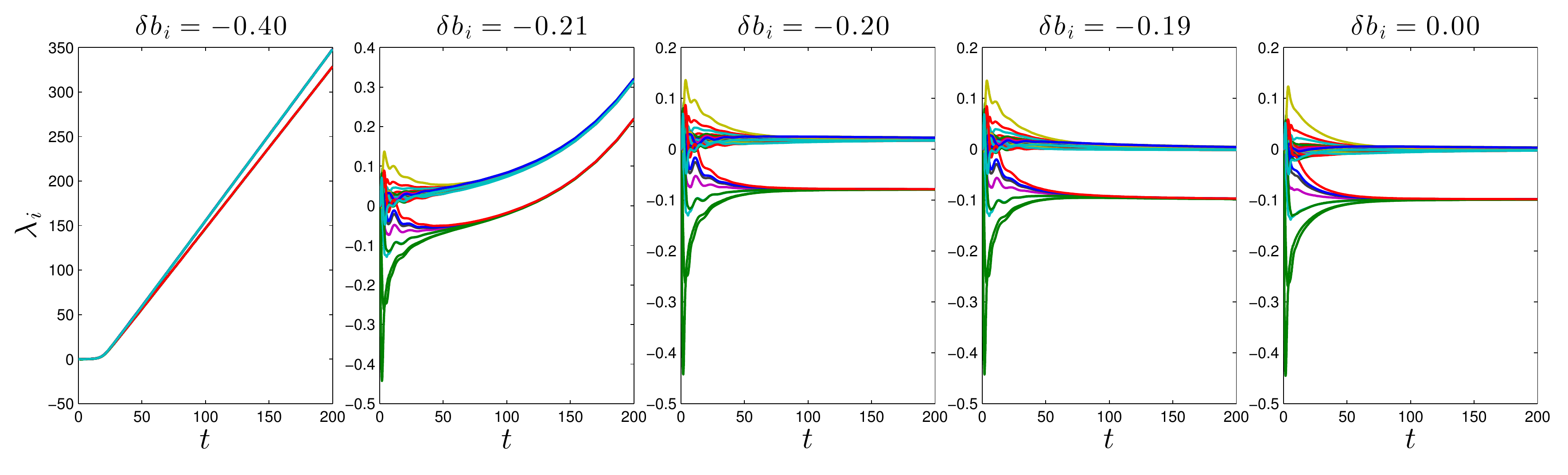}
\caption{Location Marginal Prices evolution for homogeneous perturbation $\delta a_i\in\{-0.4,-0.21,-0.2,-0.19,0.0\}$}\label{fig:perturbed-system-2}
\end{figure}

Figures \ref{fig:perturbed-system-1} and \ref{fig:perturbed-system-2} show the evolution of the frequency $\omega_i$ and LMPs $\lambda_i$
for different values of  $\delta a_i$ belonging to $\{-0.4,-0.21,-0.2,-0.19,0.0\}$. Since $D_i=0.2$ at all the buses, then $\delta a_i=-0.2$ is the threshold that makes $a_i$ go from positive to negative as $\delta a_i$ decreases. 

Even though condition \eqref{eq:condition-th16} is not satisfied for any $\delta a_i$, our simulations show that the system converges whenever $a_i\geq 0$ ($\delta a_i\geq -0.2$). The case when $\delta a_i=-0.2$ is of special interest. Here, the system converges, yet the nominal frequency is not restored. This is because the terms $\delta a_i\omega_i$ \eqref{eq:49} are equal to the terms $D_i\omega_i$ in  \eqref{eq:1a}-\eqref{eq:1b}. Thus $\dot \omega_i$ and $\dot \lambda_i$ can be made simultaneously zero with nonzero $\omega_i^*$. Fortunately, this can only happen when $a_i=0$ $\forall i$ which can be avoided since $a_i$ is a designed parameter.
\fi

%
%

\section{Concluding remarks}
\label{sec:concluding-remarks}
This paper studies the problem of restoring the power balance and operational constraints of a power network after a disturbance by dynamically adapting the loads. We show that provided communication is allowed among neighboring buses, it is possible to rebalance the power mismatch, restore the nominal frequency, and maintain inter-area flows and thermal limits. 
Our distributed solution converges for every initial condition and is robust to parameter uncertainty. Several numerical simulations verify our findings and provide new insight on the conservativeness of the theoretical sufficient condition.


\appendix

\section{Proof of Lemmas}\label{app:lemmas}

\subsection{Proof of Lemma \ref{lem:OLC-optimality}}
\begin{IEEEproof}
Assumptions \ref{as:1} and \ref{as:2} guarantee that the solution to the primal (OLC) is finite. Moreover, since by Assumption \ref{as:2} there is a feasible $d\in \text{Int}\,\D$, then the Slater condition is satisfied \cite{Boyd:2004cv} and there is zero duality gap. 

Thus, since OLC only has linear equality constraints, we can use Karush-Kuhn-Tucker (KKT) conditions~\cite{Boyd:2004cv} to characterize the primal dual optimal solution. Thus $(d^*,\omega^*,P^*,\phi^*,\sigma^*)$ is primal dual optimal if and only if we have:
\begin{enumerate}
\item[(i)] Primal and dual feasibility: \eqref{eq:olfc-b}-\eqref{eq:olfc-e} and $\rho^{+*},\rho^{-*}\geq 0$.
\item[(ii)] Stationarity:
\begin{gather*}
\frac{\partial}{\partial d} L(d^*,\omega^*,x^*,\sigma^*)=0, \quad \frac{\partial}{\partial \omega} L(d^*,\omega^*,x^*,\sigma^*)=0 \\
\text{ and }\quad \frac{\partial}{\partial x} L(d^*,\omega^*,x^*,\sigma^*)=0.
\end{gather*}
\item[(iii)] Complementary slackness: 
\begin{align*}
\rho_{ij}^{+*}(B_{ij}(v_i^*-v_j^*) - \bar P_{ij}) = 0,\qquad ij\in\E ;\\
\rho_{ij}^{-*}(\underline P_{ij}-B_{ij}(v_i^*-v_j^*)) = 0,\qquad ij\in\E.
\end{align*}
\end{enumerate}
KKT conditions (i) and (iii) are already implicit by assumptions of the lemma. 

The stationarity condition (ii) is given by
\begin{subequations}
\begin{align}
&\frac{\partial L}{\partial d_i}  (d^*,\omega^*,P^*,\phi^*,\sigma^*) = c'_i(d_i^*) -(\nu_i^*+\lambda_i^*)=0\label{eq:pLpd}\\
&\frac{\partial L}{\partial \omega_i} (d^*,\omega^*,P^*,\phi^*,\sigma^*) =  D_i(\omega_i^*-\nu_i^*)=0\label{eq:pLpomega}\\
&\frac{\partial L}{\partial P_{ij}}  (d^*,\omega^*,P^*,\phi^*,\sigma^*) = \nu_j^*-\nu_i^* = 0 \label{eq:pLpP}\\
&\frac{\partial L}{\partial \phi_i}  (d^*,\omega^*,P^*,\phi^*,\sigma^*) =\sum_{j\in\N_i} B_{ij}(\lambda_j^*-\lambda_i^*)\nonumber\\
&\qquad\quad + \sum_{e\in\E}C_{i,e}B_e(\sum_{k\in\K} \hat C_{k,e}\pi_k^* + \rho^{+*}_e-\rho^{-*}_e)=0 \label{eq:pLpphi}
\end{align}
\end{subequations}

Since $D_i>0$ equation \eqref{eq:pLpomega} implies $\nu_i^*=\omega_i^*$. Thus, \eqref{eq:pLpomega} and \eqref{eq:pLpd} amount to the first and second conditions of \eqref{eq:16}.  Furthermore, since the graph $G$ is connected then  \eqref{eq:pLpP} is equivalent to 
\[
\nu_i^*=\hat\nu \quad \forall i\in\N.
\]
which is the third condition of  \eqref{eq:16}.


Since $c_i(d_i)$ and $\frac{D_i\omega_i^2}{2}$ are strictly convex functions, it is easy to show 
that $\nu_i^*$ and $\lambda_i^*$ are unique. 
To show $\hat\nu=0$ we use (i). Adding \eqref{eq:olfc-b} over $i\in\N$ gives
\begin{align}
0&=\sum_{i\in\N} \left( P^m_i -(d_i^*+D_i\omega_i^*) -\sum_{e\in\E} C_{ie}P_{e}\right) \nonumber\\
 &=\sum_{i\in\N} \left( P^m_i -(d_i^*+D_i\omega_i^*) \right) \ - \sum_{e=ij\in\E} \left( C_{ie}P_{e}+C_{je}P_{e} \right) \nonumber\\
 &=\sum_{i\in\N} \left( P^m_i -(d_i^*+D_i\omega_i^*) \right) 
 \label{eq:sum1}
\end{align}
and similarly  \eqref{eq:olfc-c} gives
\begin{align}\label{eq:sum2}
0&=\sum_{i\in\N} P^m_i -d_i^*
\end{align}
Thus, subtracting \eqref{eq:sum1} from \eqref{eq:sum2} gives
\begin{align*}
0=\sum_{i\in\N} D_i\omega_i^*=\sum_{i\in\N}  D_i\nu_i^*
 =\hat\nu\sum_{i\in\N} D_i
\end{align*}
 and  since $D_i>0$ $\forall i\in \N$, it follows that $\hat\nu=0$. 
\end{IEEEproof}

\subsection{Proof of Lemma \ref{lem:equivalence}}
\begin{IEEEproof}
Let $(d^*\omega^*=0,\theta^*)$ be an optimal solution of OLC. Then, by letting $\phi^*=\theta^*$ and $P^*=BC^T\theta^*$, it follows that $(d^*,\omega^*=0,\phi^*,P^*)$ is a feasible solution of VF-OLC. Suppose that $(d^*,\omega^*,\phi^*,P^*)$ is not optimal with respect to VF-OLC, then the solution $(\hat d^*,\hat\omega^*,\hat\phi^*,\hat P^*)$ of VF-OLC  has strictly lower cost $\sum_{i\in\N} c_i(\hat d_i^*)+\frac{D_i\hat\omega_i^{*2}}{2} <\sum_{i\in\N} c_i(d_i^*)$. By Lemma \ref{lem:OLC-optimality} we have that $\hat\omega^*=\0$. Then, it follows that by setting $\hat \theta^*=\hat \phi^*$, $(\hat d^*,\hat\omega^*,\hat\theta^*)$ is a feasible solution of OLC with strictly lower cost than the supposedly optimal $(d^*,\omega^*,\theta^*)$. Contradiction. Therefore $(\hat d^*,\hat\omega^*,\hat\phi^*,\hat P^*)$ is an optimal solution of VF-OLC. The converse is shown analogously.
%
\end{IEEEproof}

{
\subsection{Proof of Lemma \ref{lem:strict-concavity}}
\begin{IEEEproof}
A straightforward differentiation shows that the Hessian of $\Phi_i(\nu_i,\lambda_i)$ is given by
\begin{align}
&\dddx{(\nu_i,\lambda_i)}\Phi_i(\nu_i,\lambda_i)=
\begin{blockarray}{c@{\hspace{5pt}}cc@{\hspace{5pt}}cl}
  &\matindex{\nu_i} & \matindex{\lambda_i} & &\\
\begin{block}{[c@{\hspace{5pt}}cc@{\hspace{5pt}}c]l}
& -(d_i' + D_i) & -d_i'  & & \matindex{\nu_i}\\
& -d_i'             & -d_i'  & & \matindex{\lambda_i}\\
\end{block}
\end{blockarray}\label{eq:DPhi_i}
\end{align}
where $d_i'$ is short for $d_i'(\lambda_i+\nu_i)$ and denotes de derivative of $d_i(\cdot)=c_i^{'-1}(\cdot)$ with respect to its argument.

Since $c_i$ is strictly convex $d_i'>0$. Thus, since $D_i>0$, \eqref{eq:DPhi_i} is negative definite which implies that $\Phi_i(\nu_i,\lambda_i)$ is strictly concave. Finally, it follows from \eqref{eq:L(x,sigma)} that $L(x,\sigma)$ is strictly concave 
in $(\nu,\lambda)$.
\end{IEEEproof}
}

\subsection{Proof of Lemma \ref{lem:differentiability-nu}}
\begin{IEEEproof}
We first notice that $\nu^*_i(x,y)$, $i\in\L$, depends only on $\lambda_i$ and $C_iP:= \sum_{e\in E}C_{i,e}P_e$. Which means that $\frac{\partial}{\partial v}\nu_\L^*(x,y)=0$, $\frac{\partial}{\partial \nu_\G}\nu_\L^*(x,y)=0$, $\frac{\partial}{\partial \pi}\nu_\L^*(x,y)=0$, $\frac{\partial}{\partial \rho}\nu_\L^*(x,y)=0$ and $\ddx{\lambda_\L}\nu_\L^*(x,y)$ is diagonal.

Now, by definition of $\nu^*_\L(x,y)$, for any $i\in\L$ we have
\begin{align}
&0=\frac{\partial}{\partial \nu_i}L(x,y,\nu^*_\L(x,y)) =\nonumber\\
&=\! P^m_i \!-\!(D_i\nu^*_i(x,y) \!+\! d_i(\lambda_i+\nu^*_i(x,y))) \!-\!\sum_{e\in E}C_{i,e}P_e\label{eq:stationarity}
\end{align}
Therefore, if we fix $P$ and take the total derivative of $\frac{\partial}{\partial \nu_i}L(x,y,\nu^*_\L(x,y))$ with respect to $\lambda_i$ we obtain
\begin{align}
0&=\frac{d}{d\lambda_i}\left(\frac{\partial}{\partial \nu_i}L(x,y,\nu^*_\L(x,y))\right)\\
&=-(D_i +  d_i'(\lambda_i+\nu^*_i)) \ddx{\lambda_i}\nu^*_i -  d_i'(\lambda_i+\nu^*_i)\label{eq:lem9-1}
\end{align}
where here we used $\nu^*_i$ for short of $\nu^*_i(x,y)$. 

Now since by assumption $c_i(\cdot)$ is strongly convex, i.e. $c_i''(\cdot)\geq\alpha$, $d_i'(\cdot)= \frac{1}{c_i''(\cdot)}\leq \frac{1}{\alpha}<\infty$. Thus, $(D_i +d_i')$ is finite and strictly positive, which implies that 
\[
\ddx{\lambda_i}\nu^*_i(x,y)=-\frac{ d_i'(\lambda_i+\nu^*_i(x,y))}{(\;D_i +  d_i'(\lambda_i+\nu^*_i(x,y))\;)},\quad i\in\L.
\]
Similarly, we obtain
\[
\ddx{P}\nu^*_i(x,y)=-\frac{1}{(\;D_i +  d_i'(\lambda_i+\nu^*_i(x,y))\;)}C_i,\quad i\in\L.
\]
where $C_i$ is the $i$th row of $C$.

Finally, notice that whenever $ d'_i(\lambda_i + \nu_i^*)$ exists, then $\ddx{x}\nu^*_i$ and $\ddx{y}\nu^*_i$ also exists. 
\end{IEEEproof}

\subsection{Proof of Lemma \ref{lem:second-order-derivatives}}\label{ssec:proof-second-order-derivatives}

\begin{IEEEproof}
Using the Envelope Theorem~\cite{mas1995microeconomic} in \eqref{eq:reduced-lagrangian} we have
\[
\dd{L}{x}(x,y) = \dd{L}{x}(x,y,\nu^*_\L(x,y))
\]
which implies that 
\begin{align}
&\dd{^2L}{x^2}(x,y) = \ddx{x}\left[\dd{L}{x}(x,y,\nu^*_\L(x,y))\right]\nonumber\\
&= \dd{^2L}{x^2}(x,y,\nu^*_\L(x,y)) +  \dd{^2L}{x\partial \nu_\L}(x,y,\nu^*_\L(x,y))\ddx{x}\nu^*_\L(x,y)\nonumber\\
&=  \dd{^2L}{x\partial \nu_\L}(x,y,\nu^*_\L(x,y))\ddx{x}\nu^*_\L(x,y)\label{eq:lem10-1}.
\end{align}
where the last step follows from $L(x,\sigma)$ being linear in $x$.

Now, by definition of $\nu^*_\L(x,y)$ it follows that
\begin{equation}\label{eq:dLdnuL}
\dd{L}{\nu_\L}(x,y,\nu^*_\L(x,y)) = 0.
\end{equation}

Differentiating \eqref{eq:dLdnuL} with respect to $x$ gives
\begin{align*}
0=\dd{^2L}{\nu_\L\partial x}(x,y,\nu^*_\L(x,y)) + \dd{^2L}{\nu_\L^2}(x,y,\nu^*_\L(x,y))\ddx{x}\nu^*_\L(x,y)
\end{align*}
and therefore
\begin{align}
&\dd{^2L}{x\partial \nu_\L}(x,y,\nu^*_\L(x,y))=\left[\dd{^2L}{\nu_\L\partial x}(x,y,\nu^*_\L(x,y))\right]^T\nonumber\\
&= - \ddx{x}\nu^*_\L(x,y)^T\dd{^2L}{\nu_\L^2}(x,y,\nu^*_\L(x,y)).\label{eq:lem10-2}
\end{align}
Substituting \eqref{eq:lem10-2} into \eqref{eq:lem10-1} gives
\begin{align}\label{eq:lem10-3}
\dd{^2L}{x^2}(x,y) = - \ddx{x}\nu^*_\L(x,y)^T\dd{^2L}{\nu_\L^2}(x,y,\nu^*_\L(x,y))\ddx{x}\nu^*_\L(x,y).
\end{align}

It follows from \eqref{eq:implicit-condition} and \eqref{eq:Phi_i} that
\begin{align}
\dd{^2L}{\nu_\L^2}(x,y,\nu^*_\L(x,y)) &= \dd{^2\Phi_\L}{\nu_\L^2}(\nu^*_\L(x,y),\lambda_\L)\nonumber\\
&=-(D_\L+d'_\L).\label{eq:lem10-4}
\end{align}
Therefore, substituting \eqref{eq:dnuLdx} and \eqref{eq:lem10-4} into \eqref{eq:lem10-3} gives \eqref{eq:lem10-a}.

A similar calculation using \eqref{eq:dnuLdy} gives \eqref{eq:lem10-b}.
\end{IEEEproof}

\subsection{Proof of Lemma \ref{lem:partial-derivatives}}
\begin{IEEEproof}
By definition of $g(x,y)$ we have 
\begin{align*}
\ddx{y}g(x,y)&=
\left[ \;0 \; \delta A_\G\ddx{y}\nu_\G^T \; \delta A_\L\ddx{y}\nu^*_\L(x,y)^T \; 0 \; 0\; \right]^T 
\end{align*}

Thus, using Lemma \ref{lem:differentiability-nu} we obtain $\ddx{x}g(x,y)$. A similar computation gives $\ddx{y}g(x,y)$.
\end{IEEEproof}

\bibliographystyle{IEEEtran}
\bibliography{TAC15}

\end{document}